\documentstyle[amscd,amssymb,verbatim,diagrams,12pt]{amsart}
\newcommand{\Si}{\Sigma}

\newcommand{\Conv}{{\operatorname{Conv}}}

\newcommand{\bC}{{\bf C}}

\newcommand{\ZZ}{{\cal Z}}

\newcommand{\OO}{{\cal O}}

\newcommand{\be}{{\bf e}}

\newcommand{\G}{{\Bbb G}}

\newcommand{\hra}{\hookrightarrow}

\newcommand{\CC}{{\cal C}}
\newcommand{\UU}{{\cal U}}
\newcommand{\WW}{{\cal W}}

\newcommand{\Spec}{\operatorname{Spec}}

\newcommand{\Proj}{\operatorname{Proj}}
\renewcommand{\P}{{\Bbb P}}

\newcommand{\de}{\delta}

\newcommand{\A}{{\Bbb A}}

\numberwithin{equation}{subsection}

\newtheorem{thm}{Theorem}[subsection]
\newtheorem{prop}[thm]{Proposition}
\newtheorem{conj}[thm]{Conjecture}
\newtheorem{lem}[thm]{Lemma}
\newtheorem{cor}[thm]{Corollary}
{  \theoremstyle{definition}
\newtheorem{defi}[thm]{Definition}

\newtheorem{rem}[thm]{Remark}
\newtheorem{rems}[thm]{Remarks}
}

\newcommand{\Pf}{\noindent {\it Proof}}

\newcommand{\ov}{\overline}

\renewcommand{\AA}{{\cal A}}

\newcommand{\MM}{{\cal M}}

\newcommand{\Om}{\Omega}

\newcommand{\Aut}{\operatorname{Aut}}
\newcommand{\red}{\operatorname{red}}

\renewcommand{\a}{\alpha}

\newcommand{\om}{\omega}
\newcommand{\De}{\Delta}
\newcommand{\la}{\lambda}

\newcommand{\C}{{\Bbb C}}

\newcommand{\R}{{\Bbb R}}
\newcommand{\Z}{{\Bbb Z}}
\newcommand{\Q}{{\Bbb Q}}
\newcommand{\La}{\Lambda}

\newcommand{\wt}{\widetilde}
\newcommand{\ot}{\otimes}

\newcommand{\sub}{\subset}
\newcommand{\ed}{\qed\vspace{3mm}}

\newcommand{\sslash}{\mathbin{/\mkern-6mu/}}

\title{Moduli spaces of nonspecial pointed curves of arithmetic genus $1$}
\author{Alexander Polishchuk}
\address{University of Oregon and National Research University Higher School of Economics}
\thanks{Supported in part by the NSF grant DMS-1400390 and by the Russian Academic Excellence Project `5-100'}
\begin{document}
\begin{abstract} In this paper we study the moduli stack $\UU_{1,n}^{ns}$ 
of curves of arithmetic genus $1$ with $n$ marked points,
forming a nonspecial divisor. In \cite{P-ainf-more-pts} this stack was realized as the quotient of an explicit scheme 
$\wt{\UU}_{1,n}^{ns}$, affine of finite type over $\P^{n-1}$, by the action of $\G_m^n$ . Our main result is an explicit
description of the corresponding GIT semistable loci in $\wt{\UU}_{1,n}^{ns}$. This allows us to identify some
of the GIT quotients with some of the modular compactifications of $\MM_{1,n}$ defined in \cite{Smyth-modular}
and \cite{Smyth-I}.
\end{abstract}
\maketitle

\section*{Introduction}

Recently there was a lot of interest towards birational models of the moduli spaces of pointed curves
which themselves admit a modular interpretation (see \cite{HH09}, \cite{HH13}, \cite{FS} and references therein).
A large class of examples of such birational models was constructed
by Smyth in \cite{Smyth-modular} (for arbitrary genus $g$) and \cite{Smyth-I, Smyth-II} (for $g=1$) by introducing
appropriate stability conditions for singular pointed curves. On the other hand, by studying projective geometry of curves,
or more generally, by studying algebraic structures associated with curves, one naturally obtains birational models of $M_{g,n}$
as GIT quotients. Establishing a connection between geometrically defined stabilities and the GIT stabilities is
usually quite hard (examples of such connections are surveyed in \cite{FS}).
The goal of the present paper is to describe a family of GIT quotients (with respect to a torus action) giving modular
compactifications of $M_{1,n}$, where the stability conditions admit an explicit geometric description.

Recall that in \cite{P-ainf-more-pts} we considered the moduli scheme
$\wt{\UU}_{g,n}^{ns}$, where $n\ge g$, classifying $(C,p_1,\ldots,p_n,v_1,\ldots,v_n)$, 
where $C$ is a reduced connected projective
curve of arithmetic genus $g$, $(p_\bullet)$ are distinct smooth marked points such that $H^1(C,\OO(p_1+\ldots+p_n))=0$
and $\OO_C(p_1+\ldots+p_n)$ is ample, and $v_i$ is a nonzero tangent vector at $p_i$ for each $i$.
There is a natural action of $\G_m^n$ on $\wt{\UU}_{g,n}^{ns}$, rescaling the tangent vectors, so that $(\la_1,\ldots,\la_n)$
acts by
\begin{equation}\label{action-formula-eq}
(C,p_1,\ldots,p_n,v_1,\ldots,v_n)\mapsto (C,p_1,\ldots,p_n,\la_1^{-1}v_1,\ldots,\la_n^{-1}v_n).
\end{equation}
We set $\UU_{g,n}^{ns}=\wt{\UU}_{g,n}^{ns}/\G_m^n$ (the quotient-stack). We also have
a natural $\G_m^n$-morphism
\begin{equation}\label{Gr-map}
\pi:\wt{\UU}_{g,n}^{ns}\to G(n-g,n)
\end{equation}
to the Grassmannian of $(n-g)$-dimensional subspaces in $n$-dimensional space, associating with $(C,p_\bullet,v_\bullet)$
the kernel of the surjective map 
\begin{equation}\label{H1-boundary-map-eq}
k^n\simeq H^0(C,\OO(p_1+\ldots+p_n)/\OO)\to H^1(C,\OO),
\end{equation}
where the first isomorphism is given by the tangent vectors $(v_\bullet)$.
It is shown in \cite{P-ainf-more-pts} that the morphism $\pi$ is affine of finite type (working over $\Z[1/6]$).

In the simplest case $g=0$ the moduli space $\wt{\UU}_{0,n}^{ns}$ is exactly the base of the miniversal
deformation of the coordinate cross in the $n$-space. The corresponding GIT quotients were described in \cite[Sec.\ 5]{P-ainf}:
they are the images of $\ov{M}_{0,n}$ under the maps to products of projective spaces associated with the standard
divisor classes $\psi_i$. In particular, in the main stability chamber one gets the space of so called Boggi-stable curves 
(see \cite{Boggi}, \cite[Sec.\ 4.2.1]{FS}), which gives an algebraic realization of the Kontsevich's compactification 
$K\ov{M}_{0,n}$ defined in \cite{Kon}.

In this paper we study the case $g=1$. In this case the map \eqref{Gr-map} becomes
$$\pi:\wt{\UU}_{1,n}^{ns}\to\P^{n-1}.$$
The curves corresponding to the open subscheme $V_n=\pi^{-1}(x_1\ldots x_n\neq 0)$ are well understood, since 
$V_n$ is the moduli space of {\it strongly nonspecial} curves studied in \cite{LP} (with choices of tangent vectors
at each marked point). Here ``strongly nonspecial"
means that $H^1(C,\OO(p_i))=0$ for each $i$. The curve $C$ underlying a point in $V_n$ can be either smooth, or
nodal with no disconnecting nodes (i.e., the {\it standard $m$-gon} with $1\le m\le n$), or the {\it elliptic $m$-fold curve}, with 
$1\le m\le n$.
Here the elliptic $1$-fold curve is just the cuspidal plane cubic; the elliptic $2$-fold curve is the degenerate cubic obtained as the union of a line and a conic, tangent at one point; while for $m>2$ the elliptic $m$-fold curve 
is the union of $m$ generic projective lines through a point in $\P^{m-1}$.

We prove that the complement to $V_n$ in $\wt{\UU}_{1,n}^{ns}$
is the union of the boundary divisors, consisting of curves glued transversally
from a pair of curves $C_1$ and $C_0$ of arithmetic genus $1$ and $0$, respectively (see Proposition \ref{strat-prop}).
The difference from the standard boundary divisors in $\ov{\MM}_{g,n}$ is that we allow to glue along singular points
(but not along the marked points). This allows us to decompose any curve $C$ underlying a point in $\wt{\UU}_{1,n}^{ns}$
as follows:
\begin{equation}\label{C-decomp-eq}
C=E\cup R_1\cup\ldots\cup R_k,
\end{equation}
where $E$ is the {\it minimal elliptic subcurve}, i.e., a curve of one of the types that appear over $V_n$, and
$R_i$ are (possibly singular) rational tails, i.e., connected curves of arithmetic genus $0$. Here all the gluings
are transversal (but possibly at singular points), and different $R_i$'s are attached to $E$ at different points. 
Furthermore, $\wt{\UU}_{1,n}^{ns}$ consists of all $(C,p_\bullet,v_\bullet)$, such that $C$ admits a decomposition as above,
and there is at least one marked point on each irreducible component of $C$ (see Corollary \ref{fund-dec-cor}).

The first result of our paper is a presentation of the moduli space $\wt{\UU}_{1,n}^{ns}$ as a projective scheme over an
affine scheme.

\medskip

\noindent
{\bf Theorem A} (=Thm.\ \ref{emb-thm}). {\it Let us work over $\Z[\frac{1}{6}]$. There exists a finitely generated 
$\Z[\frac{1}{6}]$-algebra $A$ and a closed
embedding $\iota:\wt{\UU}_{1,n}^{ns}\hra \P^N_A$ such that $\iota^*\OO(1)\simeq \pi^*\OO(1)$.
}

\medskip

Combining this with the fact that there are no nonconstant global $\G_m^n$-invariant functions on $\wt{\UU}_{1,n}^{ns}$
we deduce that the GIT quotients of $\wt{\UU}_{1,n}^{ns}$ associated with linearizations of
the line bundle $\pi^*\OO(1)$ are projective over $\Z[1/6]$ (see Proposition \ref{GIT-proj-prop}).

Our second 
main result is the following explicit description of the GIT stability conditions for the $\G_m^n$-action on $\wt{\UU}_{1,n}^{ns}$.
For a rational character $\chi$ of $\G_m^n$ we denote by $\pi^*\OO(1)\ot\chi$ the twist of the standard $\G_m^n$-equivariant
structure on $\pi^*\OO(1)$ by $\chi$ (so we really consider $\pi^*\OO(N)\ot\chi^N$ for sufficiently divisible $N>0$).

\medskip

\noindent
{\bf Theorem B} (=Thm.\ \ref{GIT-thm}(i)). {\it We still work over $\Z[\frac{1}{6}]$. 
Let $\chi=\sum_i a_i\be_i$ be a rational character of $\G_m^n$. For a curve
$(C,p_\bullet,v_\bullet)\in \wt{\UU}_{1,n}^{ns}$ let $E\sub C$ be the minimal elliptic subcurve, and 
let us define three subsets of $[1,n]$ as follows:
\begin{itemize}
\item
$J$ is the set of $j$ such that $p_j\not\in E$ and there are $\ge 3$ special points on the irreducible component containing $p_j$;
\item
$I$ is the set of $i$ such that $p_i\in E$;
\item
$I_0=\emptyset$ if $E$ is at most nodal; otherwise,
$I_0$ is the set of $i\in I$ such that there are $\le 2$ special points on the irreducible component containing $p_i$.
\end{itemize}
Then the point $(C,p_\bullet,v_\bullet)$ is $\pi^*\OO(1)\ot\chi$-semistable if and only if 
all $a_i\ge 0$; $a_i=0$ for $i\not\in I\cup J$; and 
$$\sum_{i\in I_0} a_i\le 1, \ \ \sum_{i\in I}a_i\ge 1.$$
}

\medskip

In particular, we see that the stability changes only when $\chi$ passes through one of the walls $a_i=0$ or
$\sum_{i\in S} a_i=1$, where $S\sub [1,n]$.
Looking at some of the resulting stability conditions we recover two different types of modular compactifications of $\MM_{1,n}$,
defined by Smyth. First, if we take $\chi=\sum_i a_i\be_i$ with each $a_i>1$ then we recover the moduli stack
$\ov{\MM}_{1,n}(\ZZ^u)$ of $\ZZ^u$-stable curves, for the extremal assignment $\ZZ^u$ of all unmarked components (see
\cite[Ex.\ 1.12]{Smyth-modular} and Proposition \ref{new-stability-prop}). In particular, we derive that the corresponding
coarse moduli space is projective (when working over $\Spec \Z[1/6]$).

Next, recall that Smyth defined in \cite{Smyth-I} the moduli stack $\ov{\MM}_{1,n}(m)$ of $n$-pointed
{\it $m$-stable} curves of arithmetic genus $1$, where $1\le m<n$. 
In terms of the decomposition \eqref{C-decomp-eq} we can characterize
these curves as follows: $(C,p_1,\ldots,p_n)$ (where $C$ is of arithmetic genus $1$
and $(p_i)$ are distinct smooth marked points) is $m$-stable if and only if
\begin{itemize}
\item $E$ is either smooth, or nodal, or the elliptic $m'$-fold curve with $m'\le m$; $R_i$ are at most nodal,
attached to $E$ so that the points of intersection are nodes on $C$;
\item $|E\cap \ov{C\setminus E}|+|E\cap\{p_1,\ldots,p_n\}|>m$;
\item $(C,p_1,\ldots,p_n)$ has no infinitesimal symmetries.
\end{itemize}
Smyth showed that $\ov{\MM}_{1,n}(m)$ is a proper irreducible Deligne-Mumford stack, with projective
coarse moduli space $\ov{M}_{1,n}(m)$ (see \cite{Smyth-I,Smyth-II}). 

We showed in \cite{P-ainf-more-pts} that for $m\ge \lfloor \frac{n}{2}\rfloor$ there exists a regular morphism
$$\ov{\MM}_{1,n}(m)\to \UU_{1,n}^{ns}$$
mapping $(C,p_\bullet)$ to $(\ov{C},p_\bullet)$,
where $\ov{C}$ is obtained by contracting the unmarked components on $C$.
We prove that this morphism factors through the semistable locus with respect to $\pi^*\OO(1)\ot\chi$ for
some rational character $\chi$ precisely for $m=n-1$, $n-2$ and $n-3$ (see Proposition \ref{Smyth-GIT-prop} and
Remark \ref{Smyth-GIT-rem}.1). Furthermore, for $m=n-1$ and $n-2$, for appropriate $\chi$,
this morphism identifies $\ov{M}_{1,n}(m)$
with the normalization of the corresponding GIT quotient (see Proposition \ref{Smyth-GIT-prop}).

The paper is organized as follows. In Section \ref{g0-sec} we discuss the genus $0$ moduli spaces $\wt{\UU}_{0,n}^{ns}$,
complementing the results of \cite[Sec.\ 5]{P-ainf}. In particular, we recall the explicit defining equations for these
affine schemes and describe the boundary divisors in $\wt{\UU}_{0,n}^{ns}$.
In Section \ref{g1-sec} we study the geometry of the genus $1$ moduli spaces $\wt{\UU}_{1,n}^{ns}$. 
After recalling the explicit affine embeddings of the open affine charts $\wt{\UU}_{1,n}(i)\sub \wt{\UU}_{1,n}^{ns}$
in Sec.\ \ref{g1-moduli-intro-sec}, we discuss in Sec.\ \ref{g1-boundary-sec}
the boundary divisors in $\wt{\UU}_{1,n}^{ns}$ (whose complement is the open subset of strongly nonspecial curves) 
and the fundamental decomposition \eqref{C-decomp-eq}. In Sec.\ \ref{global-fun-sec} we define some global functions
and global sections of $\pi^*\OO(1)$ on $\wt{\UU}_{1,n}^{ns}$ and prove some relations between them which allow us to
deduce Theorem A in Sec.\ \ref{proj-emb-sec}. Finally, in Sec.\ \ref{stability-sec} we study GIT stability conditions for
the $\G_m^n$-action on $\wt{\UU}_{1,n}^{ns}$. After some preparations we prove Theorem B in Sec.\ \ref{stability-subsec}.
We discuss the connection to Smyth's $\ZZ^u$-stabiliity in Proposition \ref{new-stability-prop} and to Smyth's $m$-stability in
Sec.\ \ref{m-stable-sec}.

\medskip

\noindent
{\it Conventions.} Starting from Section \ref{g1-sec} we work over $\Z[1/6]$. All curves (over algebraically closed fields)
are assumed to be reduced, connected and projective. By a {\it special} or a {\it distinguished} point
on a component of a curve $C$ we mean either a marked point or a singular point of $C$. 
By the {\it rational $n$-fold singularity} we mean the singularity at the origin 
of the union of $n$ coordinate lines ({\it the coordinate cross}) in the $n$-space.
By the {\it standard $n$-gon} we mean the nodal curve of genus $1$ which is the union of $n$ projective lines,
glued to form a cycle.

\section{Curves of genus $0$}\label{g0-sec}

\subsection{The stack $\wt{\UU}_{0,n}^{ns}$ and the universal affine curve over it}

Let us denote by 
$$\wt{\CC}_{g,n}^{ns}\to \wt{\UU}_{g,n}^{ns}$$ 
the universal affine curve, i.e., the complement
to the universal sections $p_1,\ldots,p_n$ in the universal curve over $\wt{\UU}_{g,n}^{ns}$.
Note that $\wt{\CC}_{g,n}^{ns}$ is the stack parametrizing  the data
$(C,p_\bullet,v_\bullet;q)$, where $(C,p_\bullet,v_\bullet)$ is in $\wt{\UU}_{g,n}^{ns}$ and
$q$ is a (possibly singular) point of $C$, different from $p_1,\ldots,p_n$.

Recall that the stack $\wt{\UU}_{0,n}^{ns}$, for $n\ge 2$, is a scheme of finite type over $\Spec(\Z)$,
given by explicit equations (see \cite[Sec.\ 5]{P-ainf}, \cite[Thm.\ 1.2.2]{P-ainf-more-pts}). 
Namely, for $n\ge 3$ the scheme $\wt{\UU}_{0,n}^{ns}$ is given by the equations
$$(\wt{\a}_{ik}-\wt{\a}_{ij})(\wt{\a}_{jk}-\wt{\a}_{ji})=(\wt{\a}_{il}-\wt{\a}_{ij})(\wt{\a}_{jl}-\wt{\a}_{ji})$$
between the variables $(\wt{\a}_{ij})_{1\le i,j\le n,i\neq j}$ defined up to translations $\wt{\a}_{ij}\mapsto \wt{\a}_{ij}+c_i$,
and the universal curve $\wt{\CC}_{0,n}^{ns}\to \wt{\UU}_{0,n}^{ns}$ is given by the equations
$$x_ix_j=\wt{\a}_{ij}x_j+\wt{\a}_{ji}x_i+c_{ij},$$
where $c_{ij}=\wt{\a}_{ik}\wt{\a}_{jk}-\wt{\a}_{ij}\wt{\a}_{jk}-\wt{\a}_{ji}\wt{\a}_{ik}$ for any $k\neq i,j$
(see \cite[Thm.\ 5.1.4]{P-ainf}, where the variables are normalized by requiring $\wt{\a}_{i,i+1}=0$).
Thus, the scheme $\wt{\UU}_{0,n}^{ns}$
is exactly the miniversal deformation space of the rational $n$-fold singularity, 
i.e., of the coordinate cross in the $n$-space (see \cite[Ex.\ 1]{Schaps}, \cite[Sec.\ 3]{Stevens}).

In the case $n=2$ the similar presentation allows to make $\wt{\a}_{12}=\wt{\a}_{21}=0$, so $c=c_{12}$ will
be the variable on $\wt{\UU}_{0,2}^{ns}\simeq \A^1$, so that we have
$$\wt{\CC}_{0,2}^{ns}\simeq \A^2$$
and the projection $\wt{\CC}_{0,2}^{ns}\to \wt{\UU}_{0,2}^{ns}$ can be identified with
the map 
$$\A^2\to \A^1: (x,y)\mapsto xy.$$

Note that for $n=1$ the universal curve $\wt{\CC}_{0,n}^{ns}$ is still a scheme.
Indeed, the stack $\wt{\CC}_{0,1}^{ns}$ parametrizes $(C,p_1,v_1;q)$, where $C\simeq \P^1$,
$p_1\neq q$ and $v_1$ is a nonzero tangent vector at $p_1$. It is easy to see that $\wt{\CC}_{0,1}^{ns}$
is isomorphic to $\Spec(\Z)$.

Sometimes it is convenient to use the following slightly different presentation of $\wt{\CC}_{0,n}^{ns}$.

\begin{lem}\label{C-0n-lem} 
For $n\ge 1$ the stack $\wt{\CC}_{0,n}^{ns}$ is isomorphic to the affine scheme 
given by the equations
\begin{equation}\label{C-0n-eq}
\a_{ik}\a_{jk}-\a_{ij}\a_{jk}-\a_{ji}\a_{ik}=0
\end{equation}
for distinct indices $i,j,k$, between the variables $(\a_{ij})$, $1\le i,j\le n$, $i\neq j$, and the doubled universal
affine curve 
$$\wt{\CC}_{0,n}^{ns}\times_{\wt{\UU}_{0,n}^{ns}}\wt{\CC}_{0,n}^{ns}\to \wt{\CC}_{0,n}^{ns}$$
is given by the equations
$$\varphi_i\varphi_j=\a_{ij}\varphi_j+\a_{ji}\varphi_i$$
in the relative $\A^n$ over $\wt{\CC}_{0,n}^{ns}$.
\end{lem}

\Pf . The proof is similar to that of \cite[Thm.\ 5.1.4]{P-ainf}. The only difference is that in our situation
we can choose uniquely $\varphi_i\in H^0(\OO(p_i))$, which vanish at the additional point $q$ and have
given polar parts (corresponding to the chosen tangent vectors). Then the equations of
$H^0(C\setminus\{p_1,\ldots,p_n\})$ become
\begin{equation}\label{genus-0-curve-eq}
\varphi_i\varphi_j=\a_{ij}\varphi_j+\a_{ji}\varphi_i,
\end{equation}
and we get the equations \eqref{C-0n-eq} using \cite[Lem.\ 5.1.3]{P-ainf}.
\ed


\begin{lem}\label{g0-sing-lem} 
We have $\wt{\UU}_{0,3}^{ns}\simeq\A^3$ with the coordinates $(\wt{\a}_{21}, \wt{\a}_{32}, \wt{\a}_{13})$, 
while $\wt{\UU}_{0,4}^{ns}$ is the affine cone over the Segre
embedding of $\P^1\times\P^3$.
If $k$ is either a field or $\Z$ then the scheme $\wt{\UU}_{0,4}^{ns}\times\Spec(k)$ is normal, while
the scheme $\wt{\UU}_{0,5}^{ns}\times\Spec(k)$ is integral.
\end{lem}

\Pf . The assertion for $n=3$ is immediate. The cases of $n=4$ and $n=5$ follow from the results of \cite{Stevens} 
about the base of the miniversal deformation of the coordinate cross.
\ed


\subsection{Irreducibility and boundary divisors}

\begin{prop}\label{genus-0-prop}
Let $k$ be either a field or $\Z$.
 For $n\ge 2$ the schemes
$\wt{\UU}_{0,n}^{ns}\times \Spec(k)$ 
and $\wt{\CC}_{0,n}^{ns}\times \Spec(k)$ 
are irreducible (of dimensions $2n-3$ and $2n-2$ when $k$ is a field).
The scheme $\wt{\UU}_{0,n}^{ns}\times \Spec(k)$ is smooth in codimension $\le 4$,
normal in codimension $\le 6$, and reduced in codimension $\le 8$.
\end{prop}
 
\Pf . The case $n=2$ was considered above, so we assume that $n\ge 3$.

The irreducibility of $\wt{\UU}_{0,n}^{ns}$ (over a field $k$ or over $\Z$)
follows from the fact that every reduced curve of arithmetic genus $0$
is smoothable (see \cite[Ex.\ 29.10.2]{Hart}).

To deduce that $\wt{\CC}_{0,n}^{ns}$ is irreducible we use the fact that the natural projection 
$$p:\wt{\CC}_{0,n}^{ns}\to \wt{\UU}_{0,n}^{ns}$$
is flat, and that the preimage of the open subset $\wt{\UU}_{0,n}^{sm}\sub \wt{\UU}_{0,n}^{ns}$ is irreducible
being isomorphic to $\MM_{0,n+1}$. 

The last assertion follows from Lemma \ref{g0-sing-lem}. Indeed, since the deformation functor of a pointed curve
is smooth over the deformations of its singular points (see e.g., \cite[Lem.\ 2.1]{Smyth-II}),
the non-smooth (resp., non-normal, resp., non-reduced) points of $\wt{\UU}_{0,n}^{ns}$ can occur only for 
curves with a rational $m$-fold singularity, where $m\ge 4$ (resp., $m\ge 5$, resp., $m\ge 6$).
It remains to use the fact that
curves with a rational $m$-fold singularity occur in codimension $\ge 2m-3$.
\ed


Recall that for every $1\le k\le n-1$ and every partition of $[1,n]$ into the disjoint union
of two subsets $[1,n]=I\sqcup J$ with $|I|=k$, we have constructed in \cite[Sec.\ 1.3]{P-ainf-more-pts} a {\it gluing morphism}
\begin{equation}\label{gluing-g0-eq}
\rho_{I,J}:\wt{\CC}_{0,k}^{ns}\times \wt{\CC}_{0,n-k}^{ns}\to \wt{\UU}_{0,n}^{ns},
\end{equation}
sending a pair of curves $C,C'$ with marked points and tangent vectors at them
and with extra points $q,q'$ (not necessarily smooth) to
the curve obtained by gluing $C$ and $C'$ transversally, identifying $q$ and $q'$, equipped with the induced
markings and tangent vectors.
Furthermore, $\rho_{I,J}$
is a closed embedding admitting a retraction (see \cite[Ex.\ 1.3.3]{P-ainf-more-pts}).
We denote by $D_{I,J}$ the image of $\rho_{I,J}$. 

\begin{cor} The subscheme $D_{I,J}$ is an irreducible
divisor in $\wt{\UU}_{0,n}^{ns}$. Furthermore, the complement to the open locus in $\wt{\UU}_{0,n}^{ns}$, corresponding 
to smooth curves, is precisely the union of all the divisors $(D_{I,J})$.
\end{cor}


\begin{rem} It was shown in \cite[Prop.\ 5.3.1]{P-ainf} 
that for a character $\chi=\sum a_i\be_i$ of $\G_m^n$, such that all $a_i>0$,
the GIT quotient $\wt{\UU}_{0,n}^{ns}\sslash_\chi \G_m^n$ can be identified with the moduli scheme $\ov{M}_{0,n}[\psi]$
of {\it $\psi$-stable} (or {\it Boggi stable}) curves, 
i.e., $n$-pointed curves of arithmetic genus $0$ without infinitesimal automorphisms, for which
there is at least one marked point on every irreducible component. The latter scheme was first constructed by Boggi in
\cite{Boggi}, who claimed that $\ov{M}_{0,n}[\psi]$ is normal. However, his proof is incorrect, and in fact, his arguments
do not even show that $\ov{M}_{0,n}[\psi]$ is reduced. So currently it is not known whether
$\ov{M}_{0,n}[\psi]$ (or equivalently, $\wt{\UU}_{0,n}^{ns}$) is reduced.
\end{rem}

\section{Curves of genus $1$}\label{g1-sec}

\subsection{The moduli spaces $\wt{\UU}^{ns}_{1,n}$}\label{g1-moduli-intro-sec}

From now on we work over $\Z[\frac{1}{6}]$.
Recall (see \cite{P-ainf-more-pts}) that we have a natural affine morphism 
$$\pi:\wt{\UU}^{ns}_{1,n}\to\P^{n-1}$$ 
sending a $k$-point $(C,p_\bullet,v_\bullet)$ to the natural functional 
$k^n\to H^1(C,\OO)\simeq k$ (see \eqref{H1-boundary-map-eq}).

We denote by $\wt{\UU}_{1,n}(i)$, $i=1,\ldots,n$, the preimages under $\pi$
of the standard affine open subsets $U_i\sub \P^{n-1}$, so that
$\wt{\UU}_{1,n}(i)$ is the open subset of $\wt{\UU}^{ns}_{1,n}$ consisting of $(C,p_1,\ldots,p_n)$
such that the map $H^0(C,\OO(p_i)/\OO)\to H^1(C,\OO)$ is an isomorphism, or equivalently,
$H^1(C,\OO(p_i))=0$. 

Consider the open subscheme
\begin{equation}\label{Vn-eq}
V_n:=\cap_{i=1}^n \wt{\UU}_{1,n}(i)\sub \wt{\UU}_{1,n}^{ns}.
\end{equation}
We will refer to curves in $V_n$ as {\it strongly nonspecial}, since to be in $V_n$ means that each point $p_i$ defines
a nonspecial divisor. The scheme $V_n$ is closely related to the moduli scheme $\wt{\UU}_{1,n}^{sns}$ studied in \cite{LP}.
Recall that $\wt{\UU}_{1,n}^{sns}$ classifies $(C,p_\bullet,\om)$, where $C$ is of arithmetic genus $1$,
$p_i$ are smooth marked points such that $H^1(C,\OO(p_i))=0$, $\OO_C(p_1+\ldots+p_n)$ is ample, 
and $\om$ is a nonzero global section of the dualizing
sheaf on $C$. The scheme $\wt{\UU}_{1,n}^{sns}$ is equipped with a natural $\G_m$-action, rescaling $\om$.
There is a canonical closed embedding
\begin{equation}\label{sns-emb-eq}
\wt{\UU}_{1,n}^{sns}\to V_n,
\end{equation}
equivariant with respect to the diagonal embedding $\G_m\to \G_m^n$, defined by choosing all tangent vectors
$v_i$ to be compatible with $\om$.
Furthermore, one has
$$V_n\simeq \wt{\UU}_{1,n}^{sns}\times_{\G_m}\G_m^n$$
(see \cite[Prop.\ 3.3.1]{P-krich}).

It follows from the results of \cite{LP} that a curve $C$ is the underlying curve of a 
strongly nonspecial curve $(C,p_\bullet)$ if and only if $C$ is either smooth, or
isomorphic to the standard $m$-gon or to the elliptic $m$-fold curve with $m\le n$ (see the proof of
\cite[Thm.\ 1.5.7]{LP}).
Furthermore, to extend such $C$ to a strongly nonspecial curve,
the (distinct smooth) marked points $p_1,\ldots,p_n$ can be chosen arbitrarily in such a way that there is at least one
on each irreducible component. This immediately leads to following assertion (which can also be checked directly).

\begin{lem}\label{sns-lem} 
Let $C$ be a curve of arithmetic genus $1$, which is either smooth, or the standard $m$-gon, or 
the elliptic $m$-fold curve. Then for every smooth point $p$ of $C$ one has $H^1(C,\OO(p))=0$.
\end{lem}

Let us recall the construction of affine embeddings of $\wt{\UU}_{1,n}(i)$ from \cite[Thm.\ 1.2.2]{P-ainf-more-pts}.
(with some simplifications due to the fact that we work over $\Z[1/6]$).
We consider the following rational functions on the universal affine curve over
$\wt{\UU}_{1,n}(i)$: $f_i\in H^0(C,\OO(2p_i))$, $h_i\in H^0(C,\OO(3p_i))$ and $h_{ij}\in H^0(C,\OO(p_i+p_j))$
for $j\neq i$. Furthermore, using the case $n=2$ (see \cite[Sec.\ 3.1]{P-krich}, \cite[Sec.\ 1.4]{P-ainf-more-pts}),
we can normalize these functions uniquely by the form of their Laurent expansions
$$f_i\equiv \frac{1}{t_i^2}+\ldots, \ h_i\equiv \frac{1}{t_i^3}+\ldots, \ h_{ij}\equiv \frac{1}{t_j}+\ldots,$$
where $(t_j)$ are formal parameters at the marked points compatible with the given tangent vectors, and 
by the condition that
the following relations hold:
\begin{equation}\label{g1-curve-eq}
\begin{array}{l}
h_i^2=f_i^3+\pi_if_i+s_i,\\
f_ih_{ij}=b_{ij}h_{ij}+a_{ij}h_i+a_{ij}e_{ij},\\
h_ih_{ij}=e_{ij}h_{ij}+a_{ij}f_i^2+a_{ij}b_{ij}f_i+a_{ij}(\pi_i+b_{ij}^2),
\end{array}
\end{equation}
where
\begin{equation}\label{s-for-eq}
s_i=e_{ij}^2-b_{ij}(\pi_i+b_{ij}^2).
\end{equation}
Note that $h_{ij}\equiv \frac{a_{ij}}{t_i}+\ldots$ at $p_i$, while
$b_{ij}=f_i(p_j)$, $e_{ij}=h_i(p_j)$.
In addition, we should have relations of the form 
\begin{equation}\label{h-ij-ij'-eq}
h_{ij}h_{ij'}=c_{j'j}(i)h_{ij}+c_{jj'}(i)h_{ij'}+a_{ij}a_{ij'}f_i+ d_{jj'}(i),
\end{equation}
for $j\neq j'$, where 
$$c_{jj'}(i)=h_{ij}(p_{j'}).$$
Using the description of $\wt{\UU}_{1,n}(i)$ in the case $n=3$ (see \cite[Prop.\ 3.2.1]{P-krich})
we deduce the equations
\begin{equation}\label{d-jj'-i-eq}
d_{jj'}(i)=a_{ij}a_{ij'}(b_{ij}+b_{ij'}).
\end{equation}
for $j\neq j'$.
Later we will calculate explicitly the restrictions of the functions $f_i$, $h_i$ and $h_{ij}$ to the minimal elliptic subcurve
$E\sub C$ in the case when $E$ is singular (see Lemma \ref{f-h-h-lem}).

Recall that the $\G_m^n$-action on $\wt{\UU}_{1,n}^{ns}$ is given by \eqref{action-formula-eq}.
The induced action of $\G_m^n$ on
$\OO(\wt{\UU}_{1,n}(i))$ is defined by $f\mapsto (\la^{-1})^*f$. We identify characters of $\G_m^n$ with the lattice
$\Z^n$ with the standard basis $\be_i$. Thus, we say that a function $f$ has $\G_m^n$-weight $\sum m_i\be_i$ if
$$f(\la^{-1}x)=\la_1^{m_1}\ldots\la_n^{m_n}\cdot f(x).$$

The proof of \cite[Thm.\ 1.2.2]{P-ainf-more-pts} together with formulas \eqref{g1-curve-eq}--\eqref{d-jj'-i-eq}
implies the following generation result.

\begin{prop}\label{generation-prop} 
Let us work over $\Z[\frac{1}{6}]$.
The ring of functions $\OO(\wt{\UU}_{1,n}(i))$ is generated
by the functions $a_{ij},b_{ij},e_{ij},\pi_i$ and $c_{jj'}(i)$. Their $\G_m^n$-weights are given by
\begin{equation}\label{U1n-weights-eq}
wt(a_{ij})=\be_j-\be_i, \ wt(b_{ij})=2\be_i, \ wt(e_{ij})=3\be_i, \ wt(\pi_i)=4\be_i, \ wt(c_{jj'}(i))=\be_j.
\end{equation}
\end{prop}

Let $(x_i)$ be the standard basis of global sections of $\OO(1)$ on $\P^{n-1}$.
It is easy to see that we have a natural identification of $\pi^*\OO(1)$ with the line bundle
$\La^*$ (dual of the Hodge bundle), with the fiber $H^1(C,\OO)$ over $(C,p_\bullet,v_\bullet)$, so that the map
\eqref{H1-boundary-map-eq} is identified with the pull back of the canonical map $\OO^n \to \OO(1)$.
We denote still by $x_i$ the pull-back of $x_i$ under $\pi$. Thus, the value of $x_i$ at $(C,p_\bullet,v_\bullet)$
is given by the image of $1/t_i\in H^0(C,\OO(p_i)/\OO)$ under the map \eqref{H1-boundary-map-eq}.
Note that over
$\wt{\UU}_{1,n}(i)$ we have the function $h_{ij}\in H^0(C,\OO(p_i+p_j))$ which projects to
$$(\frac{a_{ij}}{t_i})+(\frac{1}{t_j})\in H^0(C,\OO(p_i+p_j)/\OO)\sub H^0(C,\OO(p_1+\ldots+p_n)/\OO).$$
Hence, this element lies in the kernel of the map \eqref{H1-boundary-map-eq}, and we get the equation
\begin{equation}\label{a-x-ij-eq}
x_j=-a_{ij}x_i
\end{equation}
between sections of $\OO(1)$, valid over $\wt{\UU}_{1,n}(i)$.

\begin{lem}\label{i-relations-lem} 
One has the following relations over $\wt{\UU}_{1,n}(i)$,
where different indices are assumed to be distinct:
\begin{equation}\label{c-kj-jk-eq}
a_{ij}c_{kj}(i)=-a_{ik}c_{jk}(i);
\end{equation}
\begin{equation}\label{b-c-e-eq}
(b_{ik}-b_{ij})c_{jk}(i)=a_{ij}(e_{ij}+e_{ik});
\end{equation}
\begin{equation}\label{e-c-pi-b-eq}
(e_{ik}-e_{ij})c_{jk}(i)=a_{ij}(\pi_i+b^2_{ij}+b_{ij}b_{ik}+b_{ik}^2);
\end{equation}
\begin{equation}\label{c-quadr-eq}
c_{jk}(i)c_{j'k}(i)-c_{j'j}(i)c_{jk}(i)-c_{jj'}(i)c_{j'k}(i)=a_{ij}a_{ij'}(b_{ik}+b_{ij}+b_{ij'}).
\end{equation}
\end{lem}

\Pf . The first three relations follow from the case $n=3$ considered in \cite[Prop.\ 3.2.1]{P-krich}.
The last relation is obtained by evaluating \eqref{h-ij-ij'-eq} at $p_k$ and using \eqref{d-jj'-i-eq}.
\ed

\begin{conj} For $n\ge 3$ the relations of Lemma \ref{i-relations-lem} together with \eqref{s-for-eq}
generate the ideal of relations between our generators of
$\OO(\wt{\UU}_{1,n}(i))\ot \Z[1/6]$. 
\end{conj}

By \cite[Prop.\ 3.2.1]{P-krich}, this conjecture is true in the case $n=3$. 

\subsection{Boundary divisors and irreducibility}\label{g1-boundary-sec}

Similarly to the genus $0$ case, for every $1\le k\le n-1$ and every partition of $[1,n]$ into the disjoint union
of two subsets $[1,n]=I\sqcup J$ with $|I|=k$, we consider the {\it gluing morphism}
\begin{equation}\label{gluing-eq}
\rho^{1,0}_{I,J}:\wt{\CC}_{1,k}^{ns}\times \wt{\CC}_{0,n-k}^{ns}\to \wt{\UU}_{1,n}^{ns},
\end{equation}
which is a closed embedding (see \cite[Sec.\ 1.3]{P-ainf-more-pts}).
Let us denote by $D_{I,J}\sub\wt{\UU}_{1,n}^{ns}$ the image of $\rho^{1,0}_{I,J}$.
We say that a curve in $D_{I,J}$ is {\it glued} from the corresponding curves in $\wt{\UU}_{1,k}^{ns}$ and
$\wt{\UU}_{0,n-k}^{ns}$.

\begin{lem}\label{glued-H1-lem} 
Assume that a curve $C$ is glued from $C_1$ and $C_0$, where $C_i$ has arithmetic genus $i$.
Then for a smooth point $p\in C$, such that $p\in C_1$, 
one has $H^1(C,\OO(p))=0$ if and only if $H^1(C_1,\OO(p))=0$.
\end{lem}

\Pf . Consider the exact sequence of sheaves on $C$,
$$0\to \OO_C(p)\to  \OO_{C_1}(p)\oplus \OO_{C_0}\to \OO_q\to 0$$
where $q=C_1\cap C_2$. Since $H^1(\OO_{C_0})=0$, 
the corresponding long exact sequence gives an isomorphism
$H^1(\OO_C(p))\simeq H^1(\OO_{C_1}(p))$.
\ed

\begin{lem}\label{DIJ-lem} 
One has $D_{I,J}\sub \cup_{i\in I}\wt{\UU}_{1,n}^{ns}(i)$.
For $i_0\in I$ the intersection $D_{I,J}\cap\wt{\UU}_{1,n}(i_0)$ 
is given by the equations
\begin{equation}\label{im-rho-eq}
a_{i_0j}=0, \ b_{i_0j}=b_{i_0j'}, \ e_{i_0j}=e_{i_0j'}, \ c_{ij}(i_0)=c_{ij'}(i_0), \ c_{ji}(i_0)=0,
\end{equation}
where $i\in I\setminus\{i_0\}$, $j,j'\in J$.
\end{lem}

\Pf . The first assertion follows immediately from Lemma \ref{glued-H1-lem}.

Suppose we have a curve $(C_I,p_\bullet,v_\bullet)$ in $\wt{\UU}_{1,k}^{ns}$ with an extra point $q$ (distinct from $p_i$'s),
as well as a curve $(C_J,p_{k+1},\ldots,p_n,v_{k+1},\ldots,v_n)$ in $\wt{\UU}_{0,n-k}^{ns}$ with an extra
point $q'$ (distinct from $p_j$'s). Let $i_0\in I$ be such that $(C_I,p_\bullet,v_\bullet)\in\wt{\UU}_{1,k}^{ns}(i_0)$.
Without loss of generality we can assume that $i_0=1$, $I=\{1,\ldots,k\}$, $J=\{k+1,\ldots,n\}$.

We have the standard functions $(f_1,h_1,h_{1i})$, where $i=2,\ldots,k$, on $C_I$, as well as $\varphi_j$, for
$j=k+1,\ldots,n$, on $C_J$ satisfying the relations \eqref{g1-curve-eq}
and \eqref{genus-0-curve-eq}, where $(\varphi_j)$ are normalized by $\varphi_j(q')=0$ (see Lemma \ref{C-0n-lem}).
These functions extend to  the glued curve (by constants on the different component). 
This immediately shows that the coordinates of the glued curve in $\wt{\UU}_{1,n}^{ns}(1)$ satisfy
$$b_{1j}=f_1(q), \ e_{1j}=h_1(q), \ c_{ij}(1)=h_{1i}(q)$$
for $i\in I\setminus\{1\}$, $j\in J$.
Next, for the glued curve 
we have $h_{1j}=\varphi_j$ for $j\in J$, so that $a_{1j}=0$ and $c_{ji}(1)=0$, whereas
$c_{jj'}(1)=\a_{jj'}$.

Recall (see the proof of \cite[Prop.\ 1.3.2]{P-ainf-more-pts})
that we can recover $(C_I,p_\bullet,v_\bullet;q),(C_J,p_\bullet,v_\bullet;q')$
from the glued curve $(C,p_1,\ldots,p_n,v_1,\ldots,v_n)$ using the Proj of the appropriate Rees algebras, so that
\begin{equation}\label{O-CI-eq}
\OO(C_I\setminus\{p_1,\ldots,p_k\})=\OO(C\setminus\{p_1,\ldots,p_k\}),
\end{equation}
\begin{equation}\label{O-CJ-eq}
\OO(C_J\setminus\{p_{k+1},\ldots,p_n\})=\OO(C\setminus\{p_{k+1},\ldots,p_n\}),
\end{equation}
with $q$ obtained as the image of $p_n$ under the contraction 
$C\to C_I$, and $q'$ obtained as the image of $p_1$ under the contraction $C\to C_J$.
These constructions give a well defined map
$$r:Z_{I,J}\to\wt{\CC}_{1,k}^{ns}\times \wt{\CC}_{0,n-k}^{ns},$$
where $Z_{I,J}$ denotes the closed subscheme given by $x_{k+1}=\ldots=x_n=0$
(in fact, $r$ is exactly the retraction constructed in \cite[Prop.\ 1.3.2]{P-ainf-more-pts}). 
Restricting to the open subset $\wt{\UU}_{1,n}(1)$,
we can compute this map on the locus $a_{1j}=0$, $j=k+1,\ldots,n$.
Namely, the standard functions $(f_1,h_1,h_{1i})$ on $C_I$ correspond to the same named functions on
$C$, via the identification \eqref{O-CI-eq}, while the function $\varphi_j$, for $j\in J$, gets identified with 
$h_{1j}$ (which is regular away from $p_j$ since $a_{1j}=0$) via \eqref{O-CJ-eq}, so that
$\a_{jj'}=c_{jj'}(1)$. Now the equation 
$$\rho^{1,0}_{I,J}r(z)=z$$
on $z\in Z_{I,J}$
defining the image of $\rho^{1,0}_{I,J}$ can be easily identified with the remaining equations from \eqref{im-rho-eq}.
\ed

Recall (see \eqref{Vn-eq}) that $V_n\sub \wt{\UU}_{1,n}^{ns}$ 
denotes the open subscheme $x_1\ldots x_n\neq 0$, corresponding to strongly nonspecial curves.

\begin{prop}\label{strat-prop}
The complement to $V_n$ in $\wt{\UU}_{1,n}^{ns}$ is the union of the closed subsets $D_{I,J}$.
\end{prop}

\Pf . We have to check that any point with $x_n=0$ belongs to some $D_{I,J}$ with $n\in J$.
Applying the action of $S_n$ we can assume that $x_1\neq 0$, $c_{ni}(1)=0$ for $2\le i\le m$, and
$c_{nj}(1)\neq 0$ for $m<j<n$ (for some $m$ satisfying $1\le m\le n-1$). We claim that our point then belongs to $D_{I,J}$, where
$I=[1,m]$, $J=[m+1,n]$. Indeed, let us check the equations of Lemma \ref{DIJ-lem}.
We already know that $c_{ni}(1)=0$ for $i\in I$, $i\neq 1$.
Next, we have $a_{1n}=-x_n/x_1=0$. Then by \eqref{c-kj-jk-eq}, 
$$a_{1j}c_{nj}(1)=-a_{1n}c_{jn}(1)=0,$$
so $a_{1j}=0$ for $j\in J$. Similarly, \eqref{b-c-e-eq} and \eqref{e-c-pi-b-eq} give
$$(b_{1j}-b_{1n})c_{nj}(1)=a_{1n}(e_{1n}+e_{1j})=0,$$
$$(e_{1j}-e_{1n})c_{nj}(1)=a_{1n}(\pi_1+b_{1j}^2+b_{1j}b_{1n}+b_{1n}^2)=0,$$
so $b_{1j}=b_{1n}$ and $e_{1j}=e_{1n}$ for $j\in J$.
Next, applying \eqref{c-quadr-eq} with $(i,j,j',k)=(1,j,n,i)$, where $i\in I\setminus\{1\}$, $j\in J$, we get
$$c_{nj}(1)c_{ji}(1)=0,$$
so $c_{ji}(1)=0$. Similarly, applying \eqref{c-quadr-eq} with $(i,j,j',k)=(1,n,i,j)$ we get
$$c_{nj}(1)(c_{ij}(1)-c_{in}(1))=0,$$
so $c_{ij}(1)=c_{in}(1)$. Thus, all the equations of Lemma \ref{DIJ-lem} hold, so our point is in $D_{I,J}$.
\ed

In view of possible generalizations to higher genus, 
we will also give a purely geometric proof of Proposition \ref{strat-prop}.

\begin{lem}\label{union-curves-lem}
Let $C$ be a connected reduced projective curve over an algebraically closed field $k$. 
Let $R\sub C$ be an irreducible component of $C$ and let $C'\sub C$ be the union of the remaining irreducible
components. Further, assume that for some line bundle $L$ on $C$ one has
$h^1(C,L)=0$ and $L|_R\simeq \OO_R$. Let $C'_1,\ldots,C'_r$ be the connected components of $C'$.
Then $R\simeq \P^1$, and
\begin{equation}\label{union-curves-eq}
h^1(C,\OO)=\sum_{i=1}^r (h^1(C'_i,\OO)+\ell(\xi_i)-1),
\end{equation}
where $\xi_i$ is the scheme-theoretic intersection $C'_i\cap R$.
\end{lem}

\Pf . We have an exact sequence of coherent sheaves on $C$
\begin{equation}\label{union-curves-seq}
0\to \OO_C\to \OO_{C'}\oplus \OO_R\rTo{\de} \OO_{\xi}\to 0
\end{equation}
where $\xi=C'\cap R$.
Tensoring it with $L$ and looking at the induced long exact sequence of cohomology we deduce
the surjection $0=H^1(L)\to H^1(L|_{C'})\oplus H^1(\OO_R)$. Hence, $H^1(R,\OO)=0$, so $R\simeq \P^1$
(since it is irreducible). The second assertion is obtained by considering the long exact sequence associated with
\eqref{union-curves-seq}, and by splitting the contributions of $\OO_{C'}$ and $\OO_{\xi}$ into the direct sum
over intersections with $C'_i$.
\ed

\noindent
{\it Another proof of Proposition \ref{strat-prop}.}
We can work over an algebraically closed field.
Let $(C,p_\bullet,v_\bullet)$ be in the complement
$\wt{\UU}_{1,n}^{ns}\setminus V_n$. Without loss of generality we can assume that $H^1(C,\OO(p_n))\neq 0$.
This implies that the map $H^0(C,\OO(p_n))\to H^0(C,\OO(p_n)/\OO)$ is surjective (since it is followed in the long exact 
sequence by the map $H^1(C,\OO)\to H^1(C,\OO(p_n))$ which has to be an isomorphism). Using the vanishing 
$H^1(C,\OO(p_1+\ldots+p_n))=0$ we derive that the connecting homomorphism
$$H^0(C,\OO(p_1+\ldots+p_{n-1})/\OO)\to H^1(C,\OO)$$ 
is surjective, i.e., $H^1(C,\OO(p_1+\ldots+p_{n-1}))=0$. If $\OO_C(p_1+\ldots+p_{n-1})$ is ample then 
$(C,p_1,\ldots,p_{n-1},v_1,\ldots,v_{n-1})$ belongs to $\wt{\UU}_{1,n-1}^{ns}$. If it were in $V_{n-1}$,
then by Lemma \ref{sns-lem}, we would get $H^1(C,\OO(p_n))=0$. Hence, $(C,p_1,\ldots,p_{n-1},v_1,\ldots,v_{n-1})$ is in
the complement of $V_{n-1}$, and we can use the induction.

It remains to consider the case when there exists an irreducible component $R$ of $C$, not containing any of the points
$(p_1,\ldots,p_{n-1})$. Since there exists a marked point on every component of $C$, we see that such $R$ is unique and
$p_n\in R$. Let $C'\sub C$ be the union of all other irreducible components of $C$, and let $C'_1,\ldots,C'_r$,
$\xi_1,\ldots,\xi_r$ be as in Lemma \ref{union-curves-lem}.

Note that the summands in the equation \eqref{union-curves-eq} are non-negative, and the $i$th summand
is zero if and only if $C'_i$ has arithmetic genus $0$ and $\ell(\xi)=1$. If there exists at least one such $i$ then
$C$ is glued from $C'_i$ and the union of the remaining $C'_j$ with $R$.  
Assuming that there are no such summands, since $h^1(C,\OO)=1$, we deduce that $C'$ is connected and
either $h^1(C',\OO)=1$ and $\ell(\xi)=1$, or $h^1(C',\OO)=0$ and $\ell(\xi)=2$. In the former case $C$ is glued
from $C'$ and $R$, so assume we are in the latter case. 

If the support of $\xi$ consists of two distinct points $q_1,q_2$
then $C$ is isomorphic to the chain of projective lines with some tails of arithmetic genus $0$ glued to it transversally,
so our claim follows in this case.
Otherwise, $\xi\simeq\Spec(k[t]/(t^2))$ which is embedded
in the standard way into $R\simeq\P^1$, and also is embedded in some way into $C'$. Let $q\in C'$ be the corresponding point.
Then $q$ is a rational $m$-fold point on $C'$ (with $m\ge 1$, where $m=1$ corresponds to the smooth point).
The embedding of $\xi$ corresponds to a nonzero tangent vector at $q$. Thus, if $t_1,\ldots,t_m$ are formal parameters on
the branches of $C'$ at $q$, then the formal completion of the local ring of $C$ at $q$ is of the form
\begin{align*}
\{(f_1,\ldots,f_m,f)\in k[[t_1]]\oplus\ldots\oplus k[[t_m]]\oplus k[[t]] \ |\ & f(t)=f_1(t_1)=\ldots=f_m(t_m), \\
 & f'(t)=c_1f'_1(t_1)+\ldots+c_mf'_m(t_m)\}.
\end{align*}
Rescaling the parameters $t_i$ we can assume that each coefficient $c_i$ is either $0$ or $1$. Hence,
$C$ is the transversal union of the elliptic $m'$-fold curve (with $m'\le m+1$) with some tails of arithmetic genus zero.
Since $C$ is not in $V_n$, we deduce that it is glued.
\ed

\begin{cor}\label{fund-dec-cor} 
(i) Every curve $(C,p_\bullet)$ in $\UU_{1,n}^{ns}$ has a {\em fundamental decomposition}
$$C=E\cup R_1\cup\ldots\cup R_r,$$
where $E$ is a connected subcurve of arithmetic genus $1$, either smooth, or the standard $m$-gon, or the elliptic $m$-fold curve, and $R_i$'s are connected subcurves of arithmetic genus $0$, 
attached to $E$ transversally at distinct points ($R_i$'s can be reducible). 
Note that the intersection points $E\cap R_i$ in the fundamental decomposition are allowed to be singular both on
$E$ and on $R_i$.

\noindent
(ii) Conversely, given a curve $C$ with a decomposition as in (i), equipped with distinct smooth
marked points $p_1,\ldots,p_n$, so
that there is at least one marked point on every irreducible component, then $H^1(C,\OO(p_1+\ldots+p_n))=0$,
so $(C,p_\bullet)$ is in $\UU_{1,n}^{ns}$.
\end{cor}

\Pf . (i) This follows immediately from Proposition \ref{strat-prop}.

\noindent
(ii) Let us pick a marked point $p_i\in E$. Then by Lemma \ref{sns-lem}, we have $H^1(E,\OO(p_i))=0$.
Applying Lemma \ref{glued-H1-lem} for the gluing of $C$ from $E$ and $R_1,\ldots,R_r$ we deduce that
$H^1(C,\OO(p_i))=0$. Hence, $H^1(C,\OO(p_1+\ldots+p_n))=0$.
\ed

In the context of the above Corollary we call $E$ the {\em minimal elliptic subcurve of} $C$ and $R_j$
the {\em rational tails}.

\begin{cor}\label{H1-cor} 
For $(C,p_\bullet,v_\bullet)\in \wt{\UU}_{1,n}^{ns}$ with the minimal elliptic subcurve $E\sub C$, one has $p_i\in E$
if and only if $H^1(C,\OO(p_i))=0$ (equivalently, $x_i\neq 0$ at $(C,p_\bullet,v_\bullet)$).
\end{cor}

\Pf . If $p_i\in E$ then $H^1(C,\OO(p_i))=0$ as we have seen in the proof of Corollary \ref{fund-dec-cor}(ii).
On the other hand, if $p_i\in R_j$ then the argument similar to that of Lemma \ref{glued-H1-lem} will give
$$H^1(C,\OO(p_i))\simeq H^1(E,\OO)\neq 0.$$
\ed

\begin{cor}\label{sing-cor} 
The curves parametrized by $\wt{\UU}_{1,n}^{ns}$ have only the following types of singularities:
\begin{itemize}
\item rational $m$-fold points with $m\le n+1$;  
\item elliptic $m$-fold points with $m\le n$;
\item transversal union of an elliptic $m$-fold point and a rational $m'$-fold point, where $m+m'\le n$ (and $m'\ge 1$).
\end{itemize}
\end{cor}

\begin{rem} The fundamental decomposition of Corollary \ref{fund-dec-cor} is a generalization of a similar decomposition
for Gorenstein curves of arithmetic genus $1$ established in \cite[Lem.\ 3.1]{Smyth-I}. Note that in Gorenstein case
all the rational tails $R_i$ are nodal curves attached to the elliptic subcurve $E$ at distinct points which become nodes on $C$.
Furthermore, the elliptic component has the trivial dualizing sheaf $\om_E$ (see \cite[Lem.\ 3.3]{Smyth-I}). 
This easily implies that $C$ is Gorenstein and has $\om_C\simeq \OO_C$ if and only if $C=E$.
Thus, the locus $V_n$ of strongly nonspecial curves is characterized by the condition that $C$ is Gorenstein and has
$\om_C\simeq \OO_C$.
\end{rem}

Recall that we denote by
\begin{equation}\label{p-univ-af-curve-eq}
p:\wt{\CC}_{1,n}^{ns}\to \wt{\UU}_{1,n}^{ns}
\end{equation}
the universal affine curve over $\wt{\UU}_{1,n}^{ns}$.

\begin{prop} Let $k$ be either a field or $\Z[1/6]$.
For $n\ge 1$ the schemes $\wt{\UU}_{1,n}^{ns}\times \Spec(k)$ and $\wt{\CC}_{1,n}^{ns}\times\Spec(k)$
are irreducible (of dimension $2n$ and $2n+1$ when $k$ is a field). 
The loci $D_{I,J}\times\Spec(k)\sub \wt{\UU}_{1,n}^{ns}\times \Spec(k)$ are irreducible divisors.
The scheme $\wt{\UU}_{1,n}^{ns}\times\Spec(k)$ is regular in codimension $\le 3$ and normal in codimension $\le 4$. 
\end{prop}

\Pf . The irreducibility of $\wt{\UU}_{1,n}^{ns}\times \Spec(k)$ (where $k$ is a field or $\Z[1/6]$)
follows from the fact that all curves appearing in $\wt{\UU}_{1,n}^{ns}$ have smoothable singularities,
either by Corollary \ref{sing-cor} or by the results of \cite{Stevens-red-curve-sing}.

Since the map \eqref{p-univ-af-curve-eq} is flat, and the preimage of the smooth locus is irreducible,
we deduce that $\wt{\CC}_{1,n}^{ns}\times\Spec(k)$ is irreducible. Together with 
Proposition \ref{genus-0-prop} this implies that $D_{I,J}\times\Spec(k)$ is irreducible.

Next, we observe that 
$$\dim D_{I,J}=\dim \wt{\CC}_{1,k}^{ns} + \dim \wt{\CC}_{0,n-k}^{ns}=(2k+1)+2(n-k)-2=2n-1,$$
so each $D_{I,J}$ is a divisor in $\wt{\UU}_{1,n}^{ns}$. 

As in Proposition \ref{genus-0-prop},
the last assertion follows from Corollary \ref{sing-cor} together with some information about deformations
of the relevant singularities. Namely, using \cite[Cor.\ 2.17]{Smyth-II}, one can see that curves with
an elliptic $m$-fold singular point occur in $\wt{\UU}_{1,n}^{ns}$ in codimension $\ge m+1$. 
This implies
that the transversal unions of the elliptic $m$-fold and the rational $m'$-fold points occur in codimension $\ge m+2m'+1$
(for $m'\ge 1$). Also, it is easy to see that as in the case $g=0$ the rational $m$-fold points occur in codimension
$\ge 2m-3$.

Now we use the fact (see Lemma \ref{g0-sing-lem}) that the rational $m$-fold point 
 has smooth (resp., normal) base of the miniversal deformation
for $m\le 3$ (resp., for $m\le 4$). 
Also, by \cite[Prop.\ 1.5.12, 1.5.13]{LP}, the elliptic $m$-fold point
 has smooth (resp., normal) base of the miniversal deformation
for $m\le 5$ (resp., any $m$).
Finally, we need to know that
the miniversal deformation space is normal for the transversal union of the line and the cusp (which occurs in
codimension $4$)---this follows from \cite[Prop.\ 3.2.1]{P-krich}.
\ed

Recall that in \cite[Prop.\ 1.1.5]{LP} we constructed a natural $\G_m$-equivariant isomorphism
of the affine universal curve $\wt{\CC}_{1,n}^{sns}$ over $\wt{\UU}_{1,n}^{sns}$ with $\wt{\UU}_{1,n+1}^{sns}$, 
so that the unique $\G_m$-invariant point of $\wt{\UU}_{1,n+1}^{sns}$ corresponds to 
$(E,p_1,\ldots,p_n,q)\in \wt{\CC}_{1,n}^{sns}$, where $E$ is the elliptic $n$-fold curve, with one marked point on every
component, and $q\in E$ is the singular point.

The natural inclusion \eqref{sns-emb-eq} extends to the $\G_m$-equivariant inclusion
$$\wt{\CC}_{1,n}^{sns}\to p^{-1}(V_n)\sub \wt{\CC}_{1,n}^{ns},$$
with respect to the diagonal embedding $\G_m\to \G_m^n$.
Now the closure of every $\G_m$-orbit on $\wt{\CC}_{1,n}^{sns}\simeq\wt{\UU}_{1,n+1}^{sns}$ contains the unique
$\G_m$-invariant point. Hence, the closure of every $\G_m^n$-orbit in $p^{-1}(V_n)$ contains the $\G_m^n$-orbit
of $(E,p_\bullet,v_\bullet,q)$, where $E$ is the elliptic $n$-fold curve, $q\in E$ is its singular point.
This gives the following fact about generic $\G_m^n$-orbits in divisors $D_{I,J}$ with $|J|=2$ (which is useful for understanding
some GIT quotients of $\wt{\UU}_{1,n}^{ns}$).

\begin{lem}\label{one-tail-orbit-lem} 
Let $(C,p_\bullet,v_\bullet)\in D_{I,J}$, where $|J|=2$, be such that its fundamental decomposition has exactly one
rational tail with two marked points in $J$.
Then the closure of the $\G_m^n$-orbit of $(C,p_\bullet,v_\bullet)$ contains the curve $E\cup R$,
where $E$ is the elliptic $(n-2)$-curve, with one marked point from $I$ on each component, $R\simeq\P^1$
with two marked points in $J$, and $E\cap R$ is the singular point of $E$.
\end{lem}

\subsection{Some global functions and global sections of $\pi^*\OO(1)$}\label{global-fun-sec}

Note that by the case $n=2$ (see \cite[Lem.\ 1.4.1]{P-ainf-more-pts}),
on $\wt{\UU}_{1,n}(i)\cap \wt{\UU}_{1,n}(j)$ we have 
\begin{equation}\label{h-ji-ij-eq}
h_{ji}=a_{ji}h_{ij};
\end{equation}
\begin{equation}\label{f-ji-eq}
f_j=h_{ij}^2-a_{ij}^2f_i-a_{ij}^2b_{ij};
\end{equation}
\begin{equation}\label{h-ji-eq}
h_j=h_{ij}^3-a_{ij}^3h_i-3a_{ij}^2b_{ij}h_{ij}-2a_{ij}^3e_{ij};
\end{equation}
\begin{equation}\label{b-e-pi-ji-eq}
b_{ji}=a_{ij}^2b_{ij}, \ e_{ji}=a_{ij}^3e_{ij}, \ \pi_j=a_{ij}^4\pi_i, \ s_j=a_{ij}^6s_i.
\end{equation}

\begin{lem}
For $i\neq j$, over $\wt{\UU}_{1,n}(i)\cap \wt{\UU}_{1,n}(j)$ one has
\begin{equation}\label{a-ijk-eq}
a_{ik}=-a_{ij}a_{jk};
\end{equation}
\begin{equation}\label{c-ijk-eq}
c_{ik}(j)=a_{ji}c_{jk}(i);
\end{equation}
\begin{equation}\label{c-ijk-bis-eq}
c_{ki}(j)=a_{ik}c_{ik}(j)=-a_{jk}c_{jk}(i);
\end{equation}
\begin{equation}\label{h-jk-i-eq}
h_{jk}=h_{ik}+a_{jk}h_{ij}-c_{kj}(i),
\end{equation}
where $k\neq i,j$;
\begin{equation}\label{c-ijkm-eq}
c_{km}(j)=c_{km}(i)+a_{jk}c_{jm}(i)-c_{kj}(i),
\end{equation}
where $i,j,k,m$ are distinct.
\end{lem}

\Pf . The identities \eqref{a-ijk-eq} and \eqref{c-ijk-eq} are easy consequences of \eqref{h-ji-ij-eq}
(\eqref{a-ijk-eq} also follows from \eqref{a-x-ij-eq}). The identity \eqref{c-ijk-bis-eq} follows from
\eqref{c-kj-jk-eq} and \eqref{c-ijk-eq}:
$$c_{ki}(j)=-\frac{a_{jk}}{a_{ji}} c_{ik}(j)=-a_{jk}c_{jk}(i).$$
For \eqref{h-jk-i-eq} we observe that
the difference $h_{jk}-h_{ik}$ is regular at $p_k$ and has the expansion $\frac{a_{jk}}{t_j}+\ldots$ at $p_j$.
Hence, $m_{ijk}:=h_{jk}-h_{ik}-a_{jk}h_{ij}$ is a global section of $\OO(p_i)$, hence a constant.
Thus, we have
$$h_{jk}=h_{ik}+a_{jk}h_{ij}+m_{ijk}.$$
Multiplying both sides with $h_{ji}=a_{ji}h_{ij}$ we get
$$h_{ji}h_{jk}=a_{ji}h_{ij}h_{ik}+a_{ji}a_{jk}h_{ij}^2+a_{ji}m_{ijk}h_{ij}.$$
Using \eqref{h-ij-ij'-eq} in the right-hand side, we can rewrite this as
\begin{equation}\label{hh-ijk-eq}
h_{ji}h_{jk}=a_{ji}a_{jk}h_{ij}^2+a_{ji}[c_{kj}(i)+m_{ijk}]h_{ij}+a_{ji}c_{jk}(i)h_{ik}+a_{ik}f_i+a_{ji}d_{jk}(i),
\end{equation}
On the other hand, using \eqref{h-ij-ij'-eq} together with \eqref{h-ji-ij-eq}, \eqref{f-ji-eq} and \eqref{h-jk-i-eq} we can expand the left-hand side as
\begin{align*}
&h_{ji}h_{jk}=c_{ki}(j)h_{ji}+c_{ik}(j)h_{jk}+a_{ji}a_{jk}f_j+d_{ik}(j)=\\
&a_{ji}c_{ki}(j)h_{ij}+c_{ik}(j)[h_{ik}+a_{jk}h_{ij}+m_{ijk}]+a_{ji}a_{jk}[h_{ij}^2-a_{ij}^2f_i-a_{ij}^2b_{ij}]+d_{ik}(j).
\end{align*}
Comparing with \eqref{hh-ijk-eq} and looking at the coefficient of $h_{ij}$, we get
$$c_{ki}(j)+\frac{a_{jk}}{a_{ji}}c_{ik}(j)=c_{kj}(i)+m_{ijk}.$$
Now \eqref{c-ijk-eq} gives
$$\frac{a_{jk}}{a_{ji}}c_{ik}(j)=-c_{ki}(j).$$
Hence, from the previous equation we get
$$m_{ijk}=-c_{kj}(i).$$
Finally \eqref{c-ijkm-eq} is obtained by evaluating \eqref{h-jk-i-eq} on $p_m$.
\ed


\begin{prop}\label{global-fun-sec-prop} 
(i) There are $\G_m$-invariant
global sections $\Pi\in\pi^*\OO(-4)$ and $S\in \pi^*\OO(-6)$ on $\wt{\UU}^{ns}_{1,n}$, such that
$$\Pi|_{\wt{\UU}_{1,n}(i)}=\frac{\pi_i}{x_i^4}, \ \ S|_{\wt{\UU}_{1,n}(i)}=\frac{s_i}{x_i^6}.$$
In particular, for each $i$ we have global functions
$$\Pi_i:=x_i^4\Pi, \ \ S_i:=x_i^6S.$$

\noindent
(ii) For $i\neq j$ let us set
$$B_{ij}=\begin{cases} b_{ij} & \text{ on } \wt{\UU}_{1,n}(i);\\ 
a_{ji}^2b_{ji} & \text{ on } \wt{\UU}_{1,n}(j);\\
c_{ij}(k)^2-a_{ki}^2b_{kj}-a_{ki}^2b_{ki} & \text{ on } \wt{\UU}_{1,n}(k), k\neq i,j;
\end{cases}
$$
$$E_{ij}=\begin{cases} e_{ij} & \text{ on } \wt{\UU}_{1,n}(i);\\ 
a_{ji}^3e_{ji} & \text{ on } \wt{\UU}_{1,n}(j);\\
c_{ij}(k)^3-a_{ki}^3e_{kj}-3a_{ki}^2b_{ki}c_{ij}(k)-2a_{ki}^3e_{ki} & \text{ on } \wt{\UU}_{1,n}(k), k\neq i,j.
\end{cases}
$$
Then $B_{ij}$ and $E_{ij}$ are well defined regular functions 
on $\wt{\UU}^{ns}_{1,n}$. Their weights with respect to the $\G_m^n$-action are $2\be_i$ 
and $3\be_i$, respectively. One has the following identity on $\wt{\UU}_{1,n}^{ns}$:
\begin{equation}\label{S-Pi-eq}
S_i=E_{ij}^2-B_{ij}(\Pi_i+B_{ij}^2).
\end{equation}

\noindent
(iii) For distinct $i,j,j'$ let us set
$$C_{jj'}(i)=\begin{cases} c_{jj'}(i)\cdot x_i & \text{ on } \wt{\UU}_{1,n}(i);\\ 
-c_{ij'}(j)\cdot x_j & \text{ on } \wt{\UU}_{1,n}(j);\\
-a_{j'j}c_{ij}(j')\cdot x_{j'} & \text{ on } \wt{\UU}_{1,n}(j');\\
-a_{mi} c_{jj'}(m)\cdot x_m+a_{mj} c_{ij'}(m)\cdot x_m+a_{mi} c_{ji}(m)\cdot x_m
 & \text{ on } \wt{\UU}_{1,n}(m), \ m\neq i,j,j'.
\end{cases}$$
Then $C_{jj'}(i)$ is a well defined regular section of $\pi^*\OO(1)$ on $\wt{\UU}_{1,n}^{ns}$ of weight
$\be_i+\be_j$.
\end{prop}

\Pf . (i) We only need to check
the compatibility of the formulas giving $\Pi$ and $S$ on $\wt{\UU}_{1,n}(i)$
on the intersections $\wt{\UU}_{1,n}(i)\cap \wt{\UU}_{1,n}(j)$.
But this follows from \eqref{b-e-pi-ji-eq}.
 
\noindent 
(ii) We have to check the compatibility of our formulas for $B_{ij}$
on the intersections of the affine open charts. The compatibility on
$\wt{\UU}_{1,n}(i)\cap \wt{\UU}_{1,n}(j)$ follows from \eqref{b-e-pi-ji-eq}. Next, over the intersection
$\wt{\UU}_{1,n}(i)\cap \wt{\UU}_{1,n}(k)$, where $k\neq i,j$, we have
$$f_i=h_{ki}^2-a_{ki}^2f_k-a_{ki}^2b_{ki}$$
(see \eqref{f-ji-eq}), which gives upon evaluating at $p_j$
\begin{equation}\label{b-c-ijk-eq}
b_{ij}=c_{ij}(k)^2-a_{ki}^2b_{kj}-a_{ki}^2b_{ki},
\end{equation}
as required.
The remaining two compatibilities amount to the identity
\begin{equation}\label{c-k-k'-qu-id}
c_{ij}(k)^2-a_{ki}^2b_{kj}-a_{ki}^2b_{ki}=c_{ij}(k')^2-a_{k'i}^2b_{k'j}-a_{k'i}^2b_{k'i},
\end{equation}
that should hold over $\wt{\UU}_{1,n}(k)\cap \wt{\UU}_{1,n}(k')$, and the identity
$$a_{ji}^2b_{ji}=c_{ij}(k)^2-a_{ki}^2b_{kj}-a_{ki}^2b_{ki}$$
that should hold over $\wt{\UU}_{1,n}(j)\cap \wt{\UU}_{1,n}(k)$. The latter identity follows from \eqref{b-c-ijk-eq} with
$i$ and $j$ swapped, together with \eqref{c-ijk-bis-eq}.
To check \eqref{c-k-k'-qu-id} we apply \eqref{c-ijkm-eq} to express $c_{ij}(k')$ in terms of $c_{ij}(k)$, $c_{k'j}(k)$ and
$c_{ik'}(k)$ and \eqref{b-c-ijk-eq} to express $b_{k'j}$ and $b_{k'i}$ in terms of $c_{k'j}(k)$, $c_{k'i}(k)^2$,
$b_{kj}$, $b_{ki}$ and $b_{kk'}$. The resulting identity to check becomes
\begin{equation}\label{c-quadratic-eq}
a_{k'i}c_{ij}(k)c_{k'j}(k)-c_{ij}(k)c_{ik'}(k)-a_{k'i}c_{k'j}(k)c_{ik'}(k)=-a_{ki}^2(b_{kj}+b_{ki}+b_{kk'}),
\end{equation}
which is obtained by multiplying \eqref{c-quadr-eq}, with $(i,j,j',k)$ replaced by $(k,i,k',j)$,
by $a_{k'i}$.

Now let us check the compatibility of our formulas for $E_{ij}$.
The compatibility on $\wt{\UU}_{1,n}(i)\cap \wt{\UU}_{1,n}(j)$ follows from \eqref{b-e-pi-ji-eq}, while the compatibility
on $\wt{\UU}_{1,n}(i)\cap \wt{\UU}_{1,n}(k)$ corresponds to the identity
\begin{equation}\label{e-c-ijk-eq}
e_{ij}=c_{ij}(k)^3-a_{ki}^3e_{kj}-3a_{ki}^2b_{ki}c_{ij}(k)-2a_{ki}^3e_{ki},
\end{equation}
which is obtained by evaluating the identity
$$h_i=h_{ki}^3-a_{ki}^3h_k-3a_{ki}^2b_{ki}h_{ki}-2a_{ki}^3e_{ki}$$
(see \eqref{h-ji-eq}) on $p_j$. To check the compatibility on $\wt{\UU}_{1,n}(j)\cap \wt{\UU}_{1,n}(k)$ we evaluate the identity
$$a_{ji}^3h_j=a_{ji}^3h_{kj}^3+a_{ki}^3h_k-3a_{ji}a_{ki}^2b_{kj}h_{kj}+2a_{ki}^3e_{kj}$$
on $p_i$. Taking into account \eqref{c-ijk-bis-eq}, the needed identity then follows from
\eqref{b-c-e-eq}. Finally, to check the compatibility on $\wt{\UU}_{1,n}(k)\cap \wt{\UU}_{1,n}(k')$ we use
\eqref{c-ijkm-eq}, \eqref{b-c-ijk-eq} and \eqref{e-c-ijk-eq} 
to express $c_{ij}(k')$, $b_{k'i}$, $e_{k'i}$ and $e_{k'j}$
in terms of $c_{**}(k)$, $b_{k*}$ and $e_{k*}$. The resulting identity is equivalent to the one obtained by multiplying
\eqref{c-quadratic-eq} with $3(c_{ij}(k)+a_{k'i}c_{k'j}(k))$, and applying \eqref{b-c-e-eq} three times.

To prove \eqref{S-Pi-eq} let us consider its restrictions to the open charts. The restriction to
$\wt{\UU}_{1,n}(i)$ is simply \eqref{s-for-eq}. The restriction to $\wt{\UU}_{1,n}(j)$ follows from \eqref{s-for-eq} using
\eqref{b-e-pi-ji-eq}. Finally, the restriction to $\wt{\UU}_{1,n}(k)$, for $k\neq i,j$ gives
\begin{align*}
&a_{ki}^6s_k=[c_{ij}(k)^3-a_{ki}^3e_{kj}-3a_{ki}^2b_{ki}c_{ij}(k)-2a_{ki}^3e_{ki}]^2-\\
&[c_{ij}(k)^2-a_{ki}^2(b_{kj}+b_{ki})]^3-a_{ki}^4\pi_k(c_{ij}(k)^2-a_{ki}^2(b_{kj}+b_{ki})).
\end{align*}
This can be checked by first eliminating $s_k$ and $a_{ki}\pi_k$ using \eqref{s-for-eq} and \eqref{e-c-pi-b-eq},
and then using \eqref{b-c-e-eq} several times.

\noindent
(iii) Again we have to check that the right-hand sides are compatible on the intersections.
For the intersection $\wt{\UU}_{1,n}(i)\cap \wt{\UU}_{1,n}(j)$ this follows from
\eqref{c-ijk-eq}. Over $\wt{\UU}_{1,n}(i)\cap \wt{\UU}_{1,n}(j')$ we have
$$c_{jj'}(i)\cdot x_i=-\frac{a_{ij}}{a_{ij'}}c_{j'j}(i)\cdot x_i=-a_{j'j}c_{ij}(j')\cdot x_{j'},$$
where we used \eqref{c-kj-jk-eq} and \eqref{c-ijk-eq}. 
Next, over $\wt{\UU}_{1,n}(i)\cap\wt{\UU}_{1,n}(m)$, where $m$ is different from $i$, $j$ and $j'$, we have
by \eqref{c-ijkm-eq},
$$c_{jj'}(i)\cdot x_i=-a_{mi} c_{jj'}(m)\cdot x_m+a_{mj} c_{ij'}(m)\cdot x_m+a_{mi} c_{ji}(m)\cdot x_m.
$$
The remaining compatibilities amount to the following identities:
$$a_{j'j}c_{ij}(j')\cdot x_{j'}=c_{ij'}(j)\cdot x_j \ \text{ on } \ \wt{\UU}_{1,n}(j)\cap \wt{\UU}_{1,n}(j');$$
$$-c_{ij'}(j)\cdot x_j=-a_{mi} c_{jj'}(m)\cdot x_m+a_{mj} c_{ij'}(m)\cdot x_m+a_{mi} c_{ji}(m)\cdot x_m 
\ \text{ on } \ \wt{\UU}_{1,n}(j)\cap \wt{\UU}_{1,n}(m);$$
$$-a_{j'j}c_{ij}(j')\cdot x_{j'}=-a_{mi} c_{jj'}(m)\cdot x_m+a_{mj} c_{ij'}(m)\cdot x_m+a_{mi} c_{ji}(m)\cdot x_m 
\ \text{ on } \ \wt{\UU}_{1,n}(j')\cap \wt{\UU}_{1,n}(m);$$
\begin{align*}
&-a_{mi} c_{jj'}(m)\cdot x_m+a_{mj} c_{ij'}(m)\cdot x_m+a_{mi} c_{ji}(m)\cdot x_m =\\
&-a_{m'i} c_{jj'}(m')\cdot x_{m'}+a_{m'j} c_{ij'}(m')\cdot x_{m'}+a_{m'i} c_{ji}(m')\cdot x_{m'} 
\ \text{ on } \ \wt{\UU}_{1,n}(m)\cap \wt{\UU}_{1,n}(m').
\end{align*}
All of these identities follow easily from \eqref{c-ijk-eq}, \eqref{c-ijk-bis-eq} and \eqref{c-ijkm-eq}.
\ed

\begin{defi}
Let $A\sub \OO(\wt{\UU}_{1,n}^{ns})$ be the subring generated by all the global functions 
$(B_{ij}, E_{ij}, \Pi_i)$
defined in Proposition \ref{global-fun-sec-prop}. Note that $S_i\in A$ due to the identity \eqref{S-Pi-eq}.
\end{defi}

\begin{cor}\label{global-fun-res-cor} 
(i) The image of the restriction homomorphism $A\to \OO(\wt{\UU}_{1,n}(i))$ is the subring generated by
the following functions
\begin{itemize}
\item $b_{ij}, \ e_{ij}, \ \pi_i$ (of $\G_m^n$-weights $2\be_i$, $3\be_i$ and $4\be_i$);
\item $a_{ij}^2b_{ij}, a_{ij}^3e_{ij}, \ a_{ij}^4\pi_i$ (of $\G_m^n$-weights $2\be_j$, $3\be_j$ and $4\be_j$);
\item $c_{jj'}(i)^2-a_{ij}^2b_{ij'}, \ c_{jj'}(i)^3-a_{ij}^3e_{ij'}-3a_{ij}^2b_{ij}c_{jj'}(i)$ 
(of $\G_m^n$-weights $2\be_j$ and $3\be_j$).
\end{itemize}
Here different indices are assumed to be distinct.

\noindent
(ii) For a nonempty subset $I\sub [1,n]$ let $X_I\sub\cap_{i\in I}\wt{\UU}_{1,n}(i)$ be the closed subscheme defined
by the equations $x_j=0$ for $j\not\in I$. Then the image of the restriction homomorphism
$A\to \OO(X_I)$ is the subring generated by the images of 
\begin{itemize}
\item $b_{ij}, \ e_{ij}, \ \pi_i$ with $i\in I$;
\item $c_{jj'}(i)^2, \ c_{jj'}(i)^3$ with $i\in I$, $j,j'\not\in I$,
\end{itemize}
under the restriction homomorphisms $\OO(\wt{\UU}_{1,n}(i))\to \OO(X_I)$ (where $i\in I$).
\end{cor}

\Pf . (i) The direct use of the formulas of Proposition \ref{global-fun-sec-prop} gives a similar assertion with
the third collection of functions replaced by
$$c_{jj'}(i)^2-a_{ij}^2b_{ij}-a_{ij}^2b_{ij'}, \ c_{jj'}(i)^3-a_{ij}^3e_{ij'}-3a_{ij}^2b_{ij}c_{jj'}(i)-2a_{ij}^3e_{ij}.$$
This is equivalent to our assertion since $a_{ij}^2b_{ij}$ and $a_{ij}^3e_{ij}$ are the restrictions of
$B_{ji}$ and $E_{ji}$, respectively.

\noindent
(ii) Recall that for $i\in I$, $j\not\in I$ one has $a_{ij}=0$, and so 
$$\Pi_j|_{X_I}=B_{ji}|_{X_I}=E_{ji}|_{X_I}=0.$$
Similarly, for $j,j'\not\in I$ we can compute $B_{jj'}|_{X_I}$ and $E_{jj'}|_{X_I}$ by
first restricting them to $\wt{\UU}_{1,n}(i)$ for some $i\in I$. Thus, we get
$$B_{jj'}|_{X_I}=c_{jj'}(i)^2|_{X_I}, \ \ E_{jj'}|_{X_I}=c_{jj'}(i)^3|_{X_I}.$$
\ed

\subsection{Projective (over affine) embedding}\label{proj-emb-sec}

Recall that we denote by $A\sub \OO(\wt{\UU}_{1,n}^{ns})$ the subring generated by
the elements $(B_{ij}, E_{ij}, \Pi_i)$ defined in Proposition \ref{global-fun-sec-prop}. Thus, we can now
view $\wt{\UU}_{1,n}^{ns}$ as a scheme over $\Spec(A)$.

\begin{thm}\label{emb-thm} 
The natural morphism
\begin{equation}\label{UPA-morphism}
\wt{\UU}_{1,n}^{ns}\to \P^{n-1}_A, 
\end{equation}
given by the sections $(x_1,\ldots,x_n)$ of $\pi^*\OO(1)$, is finite.
There exists a closed embedding 
\begin{equation}\label{iota-UPA-morphism}
\iota: \wt{\UU}_{1,n}^{ns}\hra \P^N_A,
\end{equation}
where $N=n+n(n-1)(n-2)$,
such that the pull-back of $\OO_{\P^N}(1)$ is $\pi^*\OO(1)$, and the pull-backs
of the homogeneous coordinates on $\P^N$ are the sections
$(x_i)_{1\le i\le n}$ and $(C_{jj'}(i))$  (over all distinct triples $i,j,j'$) of $\pi^*\OO(1)$ .
\end{thm}

\Pf . Note that both morphisms \eqref{UPA-morphism} and \eqref{iota-UPA-morphism}
are well defined since the sections $(x_i)$ have no common zeros on $\wt{\UU}^{ns}_{1,n}$.
The fact that \eqref{UPA-morphism} is finite can be checked over each open subset $(x_i\neq 0)$ in $\P^{n-1}_A$.
Namely, we have to check that $\wt{\UU}_{1,n}(i)$ is finite over the affine space 
$\A^{n-1}_A$ over $A$ with coordinates $(a_{ij})_{j\neq i}$. But the ring $\OO(\wt{\UU}_{1,n}(i))$
is generated over $A$ by the functions $(a_{ij})$ and $(c_{jj'}(i))$, since $b_{ij}$, $e_{ij}$ and $\pi_i$ are
in the image of $A$ (see Proposition \ref{generation-prop} and Corollary \ref{global-fun-res-cor}(i)). Now the 
fact that $c_{jj'}(i)^2-a_{ij}^2b_{ij'}$ is in the image of $A$ (see Corollary \ref{global-fun-res-cor}(i))
gives the needed integral dependence for $c_{jj'}(i)$.

The fact that $\OO(\wt{\UU}_{1,n}(i))$ is generated over $A$ by $(a_{ij})$ and $(c_{jj'}(i))$
also shows that the map from $\wt{\UU}_{1,n}(i)$ to the affine chart $(x_i\neq 0)$ in $\P^N_A$,
induced by $\iota$, is a closed embedding. It follows that $\iota$ is a locally closed embedding. 
Now we observe that $\wt{\UU}_{1,n}^{ns}$ is proper over $\Spec(A)$, since it is finite over $\P^{n-1}_A$.
Hence, the image of $\iota$ is closed, and so, $\iota$ is a closed embedding.
\ed

\begin{rem} One can check (by computing restrictions to each affine chart $\wt{\UU}_{1,n}(i)$)
that the identities in the ring $\OO(\wt{\UU}_{1,n}(i))$ from Sections \ref{g1-moduli-intro-sec} and \ref{global-fun-sec}
extend to the following identities over $\wt{\UU}_{1,n}^{ns}$: 
\begin{equation}\label{C-sq-identity-eq}
C_{jj'}(i)^2=B_{jj'}x_i^2+(B_{ij}+B_{ij'})x_j^2,
\end{equation}
$$B_{ij}x_j^2=B_{ji}x_i^2, \ \ E_{ij}x_j^3=E_{ji}x_i^3, \ \ \Pi_ix_j^4=\Pi_jx_i^4,$$
$$C_{jj'}(i)=-C_{ij'}(j), \ \ C_{jj'}(i)x_{j'}=C_{ij}(j')x_j,$$
$$C_{jj'}(i)x_m=C_{jj'}(m)x_i-C_{ij'}(m)x_j-C_{ji}(m)x_i,$$
$$C_{jj'}(i)^3=E_{jj'}x_i^3-E_{ij'}x_j^3+3C_{jj'}(i)B_{ij}x_j^2-2E_{ij}x_j^3,$$
$$(B_{ik}-B_{ij})C_{jk}(i)=(E_{ij}+E_{ik})x_j, \ \ (E_{ik}-E_{ij})C_{jk}(i)=(\Pi_i+B^2_{ij}+B_{ij}B_{ik}+B_{ik}^2)x_j,$$
$$C_{jk}(i)C_{j'k}(i)-C_{j'j}(i)C_{jk}(i)-C_{jj'}(i)C_{j'k}(i)=(B_{ik}+B_{ij}+B_{ij'})x_jx_{j'}.$$
Note that the first equation \eqref{C-sq-identity-eq} gives an integral dependence for each $C_{jj'}(i)$ 
over $A[x_1,\ldots,x_n]$, which will be useful in the analysis of GIT stabilities in Section \ref{stability-sec}.
\end{rem}

Using Serre's theorems about projective schemes we immediately deduce the following Corollary from
Theorem \ref{emb-thm}.

\begin{cor}\label{global-sections-cor} 
The $\Z[1/6]$-algebra $\bigoplus_{m\ge 0} H^0(\wt{\UU}^{ns}_{1,n},\pi^*\OO(m))$ is finitely generated.
Let $\AA\sub\bigoplus_{m\ge 0} H^0(\wt{\UU}^{ns}_{1,n},\pi^*\OO(m))$
be the graded $A$-subalgebra generated by $(x_i)$ and $(C_{jj'}(i))$. Then
one has $H^0(\wt{\UU}^{ns}_{1,n},\pi^*\OO(m))=\AA_m$ for $m\gg 0$.
\end{cor}

\begin{cor}\label{finite-cor}
 The morphism
$\wt{\UU}_{1,n}^{ns}\to \Spec(A)$ is surjective, projective, and finite over 
the union of the distinguished open affine subsets $\cup_i (D(\Pi_i)\cup D(S_i))\sub\Spec(A)$.
\end{cor}

\Pf . By Theorem \ref{emb-thm}, the map $\wt{\UU}_{1,n}^{ns}\to \Spec(A)$ is projective.
Hence, it is surjective (since it is dominant). Next, since $\Pi_i=x_i^4\Pi$, we deduce 
that the preimage of $D(\Pi_i)$ in $\wt{\UU}_{1,n}^{ns}$ is the 
open subset given by the nonvanishing of $x_i$ and of $\Pi$. Hence, the preimage of $D(\Pi_i)$ is affine.
Similarly, the preimage of $D(S_i)$ is affine. Thus, over the union of these open sets our morphism is both projective
and affine, hence it is finite.
\ed

\section{Stability conditions for genus $1$ curves}\label{stability-sec}

\subsection{Explicit form of the canonical generators}

In this section we continue to work over $\Z[\frac{1}{6}]$.
Recall (see Sec.\ \ref{g1-moduli-intro-sec})
that for each curve in $(C,p_\bullet,v_\bullet)$ in $\wt{\UU}_{1,n}(i)$
we have canonical generators $(f_i,h_i,h_{ij})$ in the algebra
$\OO(C\setminus\{p_1,\ldots,p_n\})$.
The following computation of the restrictions of these generators to the minimal elliptic subcurve
will play a crucial role in the analysis of stability conditions on $\wt{\UU}_{1,n}^{ns}$.

\begin{lem}\label{f-h-h-lem}
For $(C,p_\bullet,v_\bullet)\in\wt{\UU}_{1,n}^{ns}$ let $E\sub C$ be the minimal elliptic subcurve, and let $E=\cup_{i=1}^m E_i$
be its decomposition into irreducible components.

\noindent
(i) Assume that $E$ is the elliptic $m$-fold curve with the singular point $q$, and assume that $p_i\in E_1\sub E$,
so that $(C,p_\bullet,v_\bullet)\in\wt{\UU}_{1,n}^{ns}(i)$ (see Cor.\ \ref{H1-cor}). Then
we have
$$f_i|_{E_k}=\begin{cases} u_i^2, & k=1, \\ 0, & k\neq 1,\end{cases}$$
$$h_i|_{E_k}=\begin{cases} u_i^3, & k=1, \\ 0, & k\neq 1,\end{cases}$$
for a coordinate $u_i$ on the normalization of $E_1$ such that $u_i=0$ at $q$,
$u_i=\infty$ at $p_i$, and $1/u_i$ is compatible with $v_i$. If $p_j\in E_1$ corresponds to $u_i=\la$ then
$$h_{ij}|_{E_k}=\begin{cases} a_{ij}\frac{u_i^2-\la u_i+\la^2}{u_i-\la}, & k=1, \\ -\la a_{ij}, & k\neq 1.\end{cases}$$
In the case when $p_j\in E_k$, where $k\neq 1$, we have
$$h_{ij}|_{E_l}=\begin{cases} a_{ij}u_i, & l=1, \\ u_j, & l=k, \\ 0, & l\neq 1,k,\end{cases}$$
where $u_j$ is a coordinate on $E_k$, such that $u_j=0$ at $q$, $u_j=\infty$ at $p_j$, and $1/u_j$ is compatible with $v_j$.
Note that in both cases $a_{ij}\neq 0$.

\noindent
(ii) Assume that $E$ is the standard $m$-gon, with the components $E_i$ ordered cyclically, 
and assume that $p_i\in E_1$. Then 
$$f_i|_{E_k}=\begin{cases} u_i^2-\frac{2}{3}, & k=1, \\ \frac{1}{3}, & k\neq 1,\end{cases}$$
$$h_i|_{E_k}=\begin{cases} u_i(u_i^2-1), & k=1, \\ 0, & k\neq 1,\end{cases}$$
where $u_i$ is a coordinate on the normalization of $E_1$ such that $u_i=\infty$ at $p_i$, $1/u_i$ is compatible with $v_i$, 
and $u_i=\pm 1$
at the intersections of $E_1$ with $E_2$ and $E_n$ (or at the preimages of the node, if $m=1$). 
If $p_j\in E_1$ corresponds to $u_i=\la$ then we have
$$h_{ij}|_{E_k}=\begin{cases} a_{ij}\frac{u_i^2-\la u_i+\la^2-1}{u_i-\la}, & k=1, \\ -a_{ij}\la, & k\neq 1.\end{cases}$$
In the case when $p_j\in E_k$, where $k\neq 1$, we have 
$$h_{ij}|_{E_l}=\begin{cases} a_{ij}u_i, & l=1, \\ u_j, & l=k, \\ \pm a_{ij}, & l\neq 1,k,\end{cases}$$
where $u_j$ is a coordinate on  $E_j$ such that $u_j=\infty$ at $p_j$ and $1/u_j$ is compatible with $v_j$.
Here again $a_{ij}\neq 0$ in both cases.

\noindent
(iii) Assume that $p_i\in E$ and $p_j\not\in E$.  Let $C'\simeq\P^1$ be the irreducible component
of $C$ containing $p_j$, and let $q\in C'$ be the point connecting $C'$ to $E$ or to the next component in
the rational chain linking $C'$ with $E$. 
Also, let $C''$ be the connected component of $\ov{C\setminus C'}$, not containing $E$. Then
$h_{ij}|_{C\setminus(C'\cup C'')}=0$
and $h_{ij}$ does not vanish on $C'\setminus\{q,p_j\}$.
\end{lem}

\Pf . (i) Since the parameter $t_i=1/u_i$ is compatible with $v_i$, we should have 
$$f_i|_{E_1}=u_i^2+a, \ \ h_i|_{E_1}=u_i^3+bu_i^2+c$$ 
(the absence of the terms with $u_i$ is dictated by the condition that $f_i$ and $h_i$ are constant on other components of $E$,
together with the fact that $E$ has an elliptic $m$-fold singularity at $q$). Now it is easy to see that the equation of the form 
$$h_i^2=f_i^3+\pi_if_i+s_i,$$
holds only if $a=b=c=0$. Alternatively, one can observe that
$f_i$ and $h_i$ are the pull-backs of similar functions on the cuspidal curve under the natural map
contracting all components except $E_1$. 

Recall that we normalize $h_{ij}\in H^0(C,\OO(p_i+p_j)$ 
by the requirement that $h_{ij}\equiv \frac{a_{ij}}{t_i}+\ldots$ near $p_i$ (where $t_i=1/u_i$) and that
the equation
\begin{equation}\label{h1j-normalization-eq}
f_ih_{ij}=b_{ij}h_{ij}+a_{ij}h_i+a_{ij}e_{ij}
\end{equation}
holds (i.e., there is no term with $f_i$ in the right-hand side). 

Assume first that $p_j\in E_k$, where $k\neq 1$. Then $b_{ij}=f_i(p_j)=0$ and $e_{ij}=h_i(p_j)=0$.
Hence, restricting \eqref{h1j-normalization-eq} to $E_1$ we get
$$h_{ij}|_{E_1}u_i^2=a_{ij}u_i^3,$$
which gives $h_{ij}|_{E_1}=a_{ij}u_i$. In particular, $h_{ij}(q)=0$, so
$h_{ij}$ vanishes on $E\setminus (E_1\cup E_k)$.
Now the formula for $h_{ij}|_{E_k}$ is dictated by the fact that it is a section of $\OO_{E_k}(p_j)$
that vanishes at $q$ and has an expansion
$\frac{1}{t_j}+\ldots$ near $p_j$, where $t_j=1/u_j$.

In the case when the point $p_j$ is in $E_1$ and corresponds to $u_i=\la$ we get from \eqref{h1j-normalization-eq} the equation
$$h_{ij}|_{E_1}\cdot u_i^2=\la^2h_{ij}|_{E_1}+a_{ij}u_i^3+a_{ij}\la^3,$$
i.e.,
$$h_{ij}|_{E_1}=a_{ij}\frac{u_i^3+\la^3}{u_i^2-\la^2}=a_{ij}\frac{u_i^2-\la u_i+\la^2}{u_i-\la}.$$
In particular, $h_{ij}(q)=-\la a_{ij}$, which implies that the restriction of $h_{ij}$ to other components is given by the
same constant.

\noindent
(ii) Contracting all the components of $E$ except for $E_1$ gives the nodal cubic curve $y^2=x^3+x^2$, with
$x=u_i^2-1$, $y=u_i(u_i^2-1)$, for some coordinate $u_i$ on $E_1$ (or the normalization of $E_1$, if $n=1$)
such that $u_i=\infty$ at $p_i$. To get rid of the quadratic term in $x$ we set
$$f_i|_{E_1}=x+\frac{1}{3}=u_i^2-\frac{2}{3},$$
so that $f_i$ and $h_i|_{E_1}=y$ satisfy the equation of the required form
$$h_i^2=f_i^3-\frac{1}{3}f_i+\frac{2}{27}.$$
We then extend $f_i$ and $h_i$ to the entire $C$ (by constants on other components).
Note that the nodes of $E$ connecting $E_1$ to other components (or to itself if $n=1$) correspond to $x=y=0$,
i.e., $u_i=\pm 1$. Thus, we have $f_i=1/3$ and $h_i=0$ at these points.

Now let us consider the functions $h_{ij}|_E$.
In the case when $p_j\in E_k$ with $k\neq 1$, we have $h_{ij}|_{E_1}=au_i+b$ for some constants $a\neq 0$ and $b$.
Note that in this case $b_{ij}=f_i(p_j)=1/3$, so $f_i-f_i(p_j)=x=u_i^2-1$.
The normalization that $f_ih_{ij}$ has no coefficient with $f_i$ is satisfied for $b=0$:
$$(f_i-f_i(p_j))h_{ij}|_{E_1}=(u_i^2-1)au_i=ah_i.$$
Thus, we have $h_{ij}|_{E_1}=a_{ij}u_i$. In particular, the value of $h_{ij}$ at the nodes connecting $E_1$
to $E_2$ and $E_n$ is $\pm a_{ij}$. These are also the values of the restrictions of $h_{ij}$ to the components $E_i$
different from $E_1$ and $E_k$. The restriction $h_{ij}|_{E_k}$ can be taken as a coordinate on $E_k$.

In the case when $p_j\in E_1$ and $u_i=\la$ at $p_j$,
we should have $h_{ij}|_{E_1}=(au_i^2+bu_i+c)/(u_i-\la)$, with $a\neq 0$.
In this case $f_i(p_j)=\la^2-2/3$, so 
$$(f_i-f_i(p_j))h_{ij}|_{E_1}=(u_i^2-\la^2)h_{ij}|_{E_1}=(u_i+\la)(au_i^2+bu_i+c).$$
The required normalization is that this is a linear multiple of $h_i=u_i(u_i^2-1)$ plus a constant,
so we derive that
$$h_{ij}|_{E_1}=a\frac{u_i^2-\la u_i+\la^2-1}{u_i-\la},$$
and the expansion $h_{ij}\equiv a_{ij}{t_i}+\ldots$ near $p_i$ (where $t_i=1/u_i$) implies that $a=a_{ij}$.
Substituting $u_i=\pm 1$ we get the value $-a_{ij}\la$ at the nodes, and hence at the remaining components of $E$.

\noindent
(iii) Note that $a_{ij}=0$ since $p_j\not\in E$ (see Corollary \ref{H1-cor}). Note also that 
$h_{ij}|_E$ is a section of $H^0(E,\OO(p_i))$, hence a constant (see Lemma \ref{sns-lem}). Therefore,
$h_{ij}|_E=h_{ij}(q)$. Now restricting the equation
$$f_ih_{ij}=b_{ij}h_{ij}+a_{ij}h_i+a_{ij}e_{ij}=b_{ij}h_{ij}$$
 to $E$ we get that
$h_{ij}(q)f_i|_E$ is constant, which implies that $h_{ij}|_E=h_{ij}(q)=0$. Hence, $h_{ij}$
vanishes on the rational chain linking $C'$ with $E$ and on every rational tail of $C$ not containing $p_j$.
Since $h_{ij}|_{C'}$ is a section of $\OO(1)$ with the pole
at $p_j$ and the zero at $q$, it does not vanish on $C'\setminus\{q,p_j\}$.
\ed

\subsection{Analysis of stability}\label{stability-subsec}


We refer the reader to \cite[Sec.\ 11]{Alper} for the general theory of GIT stability and GIT quotients over an arbitrary 
base scheme (recall that in our case the base is $\Spec(\Z[1/6])$).

We are going to study the GIT stabilities on $\wt{\UU}_{1,n}^{ns}$ with respect to certain $\G_m^n$-equivariant ample
line bundles. The notion of GIT stability does 
not change if we replace a $\G_m^n$-equivariant line bundle by its positive power, so
it is convenient to work with rational characters of $\G_m^n$, $\chi=\sum a_i\be_i\in\Q^n$. For such $\chi$ the notion
of $\pi^*\OO(1)\ot\chi$-(semi-)stability is defined to be that of $\OO(N)\ot (N\chi)$-(semi-)stability, where $N>0$ is such that
$N\chi\in \Z^n$. We denote by 
$\wt{\UU}^{ns}_{1,n}\sslash_\chi \G_m^n$ the GIT quotient corresponding to $\pi^*\OO(1)\ot\chi$-semistability.

Note that the action of $\G_m^n$ on the global sections of equivariant line bundles is given by 
$s\mapsto (\la^{-1})^*s$. In particular, $\G_m^n$-invariant sections of $\OO(N)\ot (N\chi)$ are identified with
sections of $\OO(N)$ of $\G_m^n$-weight $N\chi$.

Before analyzing in detail the stability conditions, we make the following general observation.

\begin{prop}\label{GIT-proj-prop} All the GIT quotients $\wt{\UU}^{ns}_{1,n}\sslash_\chi \G_m^n$ are projective
over $\Z[1/6]$. 
\end{prop}

\Pf . Indeed, since $\wt{\UU}^{ns}_{1,n}$ is projective over $\Spec(A)$ (by Theorem \ref{emb-thm}), the GIT quotient
$\wt{\UU}^{ns}_{1,n}\sslash_\chi \G_m^n$
can be identified with
$$\Proj\left(\bigoplus_{m\ge 0} H^0(\wt{\UU}^{ns}_{1,n},\pi^*\OO(m))_{\chi^m}\right),$$
which is projective over $\Spec\left(H^0(\wt{\UU}^{ns}_{1,n},\OO)^{\G_m^n}\right)$.
It remains to prove that for each $i$ we have
$$H^0(\wt{\UU}^{ns}_{1,n}(i),\OO)^{\G_m^n}=\Z[1/6].$$
This is proved similarly to \cite[Lem.\ 2.2.3(ii)]{P-krich}. Let us consider the functional $\ell$ on $\R^n$ given by
$\ell(\be_i)=1$, $\ell(\be_j)=2$ for $j\neq i$. Then our assertion follows from
the fact that all generators of $H^0(\wt{\UU}^{ns}_{1,n}(i),\OO)$ from Proposition \ref{generation-prop} have $\G_m^n$-weights in
the half-space $\ell>0$.
\ed

For each set of vectors $\Om\sub\R^n$ we denote by $\bC(\Om)$ and by $\Conv(\Om)$
the closed cone generated by $\Om$ and the convex hull of $\Om$, respectively.


\begin{defi}\label{conv-Cp-def}
For each $(C,p_\bullet,v_\bullet)$ in $\wt{\UU}^{ns}_{1,n}$ 
we denote by $\Conv(C,p_\bullet)$ the convex hull in $\R^n$ of the vectors 
$wt/N$, for all $N>0$, all $\G_m^n$-weights $wt$, and all global sections $s$ of $\pi^*\OO(N)$, 
such that $s$ has weight $wt$ with respect to the $\G_m^n$-action and $s$ is nonzero at $(C,p_\bullet,v_\bullet)$.
\end{defi}

Let $\chi=a_1\be_1+\ldots+a_n\be_n$ be a rational character of $\G_m^n$.
Since the morphism $\pi$ is affine, it follows easily from the definition 
that a point $(C,p_1,\ldots,p_n)$ is $\pi^*\OO(1)\ot\chi$-semistable if and only if
$\chi\in \Conv(C,p_\bullet)$ (we use the fact that an open subset in $\wt{\UU}^{ns}_{1,n}$, given by the nonvanishing
of a section of $\pi^*\OO(N)$, is affine). Thus, to analyze the stability conditions we have to compute
$\Conv(C,p_\bullet)$.

\begin{defi}
For $(C,p_\bullet,v_\bullet)$ in $\wt{\UU}^{ns}_{1,n}$ let 
$\Om_0(C,p_\bullet)$ denote the set of weights of the generators $(B_{ij}, E_{ij}, \Pi_i)$ of $A$, 
that are not zero at $(C,p_\bullet)$.
Similarly, let $\Om_1(C,p_\bullet)$ denote the set of $\be_i$ such that $x_i\neq 0$ at $(C,p_\bullet)$.
\end{defi}

\begin{lem}\label{cones-lem1} 
One has 
\begin{equation}\label{Conv-eq}
\Conv(C,p_\bullet)=\bC(\Om_0(C,p_\bullet))+\Conv(\Om_1(C,p_\bullet)).
\end{equation}
\end{lem}

\Pf . By Corollary \ref{global-sections-cor}, in the definition of
$\Conv(C,p_\bullet)$ (see Def.\ \ref{conv-Cp-def}) it is enough to consider $s\in\AA_N$.
Furthermore, it is enough to consider $s$ which are monomials in $(x_i)$ and $(C_{jj'}(i))$ with coefficients in $A$. The identity 
\eqref{C-sq-identity-eq} shows that we get the same set by considering only $s$ of the form $a\cdot M$, where
$M$ is a monomial of positive degree in $(x_i)$ and $a$ is a monomial in generators of $A$. We have 
$$\frac{wt(a\cdot M)}{\deg(M)}=\frac{wt(a)}{\deg(M)}+\frac{wt(M)}{\deg(M)},$$
which shows that the left-hand side of \eqref{Conv-eq} is contained in the right-hand side.
Now assume that $a\in A$ and $x_i$ do not vanish at $(C,p_\bullet)$. Then $\be_i\in \Conv(C,p_\bullet)$, and
$$wt(a^m\cdot x_i)=m\cdot wt(a)+\be_i\in \Conv(C,p_\bullet)$$
for any $m\ge 0$. Passing to the convex hull, we derive that the right-hand side of \eqref{Conv-eq} is contained in
the left-hand side.
\ed

Let $C=E\cup R_1\cup\ldots\cup R_k$ be the fundamental decomposition of a curve
$(C,p_\bullet)$ in $\UU_{1,n}^{ns}$ (see Corollary \ref{fund-dec-cor}). 

\begin{defi}\label{attached-defi}
Given an irreducible component $E'$ of $E$,
let us say that a marked point $p_j$ is {\it attached to } $E'$ if either $p_j\in E'$, or
there exists a rational tail $R_j$ such that $p_j\in R_j$ and $R_j$ is attached to $E'$ in such a way
that $R_j\cap E'$ is a smooth point of $E$.
\end{defi}

\begin{lem}\label{cones-lem2} 
Let $(C,p_\bullet,v_\bullet)$ be in $\wt{\UU}_{1,n}^{ns}$, and let $E\sub C$ be the minimal elliptic subcurve.

\noindent
(i) $E$ has at most nodal singularities if and only if either $\Pi\neq 0$ or $S\neq 0$ at $(C,p_\bullet,v_\bullet)$.

\noindent
(ii) One has $\Om_1(C,p_\bullet)=\{\be_i \ |\ p_i\in E\}$.

\noindent
(iii) Let $\Om'_0(C,p_\bullet)$ be the set of $\be_j$ such that $p_j\not\in E$, and there are at least $3$ distinguished
points on the irreducible component of $C$ containing $p_j$. 

If $E$ is at most nodal then the cone $\bC(\Om_0(C,p_\bullet))$ is spanned by $\Om'_0(C,p_\bullet)$ and by
all $\be_i$ such that $p_i\in E$.

If $E$ is the elliptic $m$-fold curve then $\bC(\Om_0(C,p_\bullet))$ is spanned by $\Om'_0(C,p_\bullet)$ and
by $\be_i$, such that $p_i\in E'\sub E$, where $E'$ is an irreducible component of $E$, and there exists another marked point attached to $E'$.
\end{lem}

\Pf . (i) This immediately follows from Lemma \ref{f-h-h-lem}.

\noindent
(ii) Recall that the condition $p_i\in E$ is equivalent to the condition $h^1(C,\OO(p_i))=0$, or $x_i\neq 0$ (see Corollary \ref{H1-cor}). Hence, $p_i\in E$ precisely when $(C,p_\bullet)$ belongs to the open subset $x_i\neq 0$.

\noindent
(iii) Let $I$ be the set of $i$ such that $x_i\neq 0$ at $(C,p_\bullet)$. 
It follows from Corollary \ref{global-fun-res-cor}(ii)
that the cone $\bC(\Om_0(C,p_\bullet))$ is spanned by 
\begin{itemize}
\item $\be_i$, such that $i\in I$, and one of the functions $b_{ij}, e_{ij}, \pi_i$ is not zero at $(C,p_\bullet)$;
\item $\be_j$, such that $j\not\in I$, and one of the functions $c_{jk}(i)$ (with $i\in I$, $k\not\in I$) is not zero  at $(C,p_\bullet)$.
\end{itemize}

If $E$ is smooth or nodal, and $i\in I$ then $p_i\in E$, so we have either $\pi_i\neq 0$ or $s_i\neq 0$, which implies that 
$\be_i$ is in $\bC(\Om_0(C,p_\bullet))$.

Now assume that $E$ is the elliptic $m$-fold curve and $i\in I$, and let $E'\sub E$ be the irreducible component
containing $p_i$. Then $\pi_i=0$, and
by Lemma \ref{f-h-h-lem}, the functions $f_i$ and $h_i$ vanish on irreducible components of $E$ other than $E'$,
and are nonzero on all smooth points of $E$ that belong to $E'$. 
This implies that $b_{ij}=e_{ij}=0$ unless $p_j$ is attached to $E'$.

In either case, if $p_i\in E$ and $p_j\not\in E$ then we claim that $c_{jk}(i)=h_{ij}(p_k)\neq 0$ for some other marked point 
$p_k\not\in E$
if and only if there are at least $3$ distinguished points on $C'$, the irreducible component containing $p_j$. 
Indeed, let $q\in C'$ be the point connecting $C'$ to $E$ or to the next component in
the chain linking $C'$ with $E$. By Lemma \ref{f-h-h-lem}(iii), we have $h_{ij}(q)=0$.
This shows that $c_{jk}(i)$ vanishes in the case when $p_j$ and $q$ are the only distinguished points of $C_j$.
On the other hand, since $h_{ij}$ is nonzero on $C'\setminus\{p_j,q\}$ (see Lemma \ref{f-h-h-lem}(iii)),
if $q'\in C'\setminus\{p_j,q\}$ is another distinguished point then $h_{ij}(q')\neq 0$.
Hence, in this case we can find a marked point $p_k\not\in E$ with
$h_{ij}(p_k)=h_{ij}(q')\neq 0$:  if $q'$ is a marked point then we take $p_k=q'$, otherwise, we take $p_k$ to be
any marked point on the component attached to $q'$.
\ed

\begin{lem}\label{Aut-curve-lem} 
Let $(C,p_\bullet,v_\bullet)$ be in $\wt{\UU}_{1,n}^{ns}(k)$, where $k$ is an algebraically closed field of characteristic
$\neq 2,3$, and let $E\sub C$ be the minimal elliptic subcurve.
Then the following conditions are equivalent:

\noindent
(a) the stabilizer subgroup scheme in $\G_m^n$ of the point $(C,p_\bullet,v_\bullet)\in\wt{\UU}_{1,n}^{ns}$ is finite and reduced;

\noindent
(b) the pointed curve $(C,p_\bullet)$ has no infinitesimal automorphisms;

\noindent
(c) every component of $C$ not contained in $E$ has $\ge 3$ distinguished points, and if $E$ is non-nodal then
there exists at least one irreducible component of $E$ with $\ge 3$ distinguished points.
\end{lem}

\Pf . It is easy to see that 
the group scheme $\Aut(C,p_\bullet)$ is isomorphic to the stabilizer subgroup of $(C,p_\bullet,v_\bullet)$ in $\G_m^n$.
This implies the equivalence of (a) and (b). The equivalence with (c) follows from the fact that an elliptic $m$-fold curve $E$,
equipped with one marked point on each component has the group $\G_m$ as automorphisms, but
there are no infinitesimal automorphisms preserving in addition one more smooth point of $E$ (see \cite[Prop.\ 1.5.4]{LP};
note that this uses the assumption on the characteristic of $k$).
\ed

For each rational character $\chi=\sum a_i\be_i$ of $\G_m$ let us denote by
$\wt{\UU}_{1,n}^{ns}(\chi)\sub \wt{\UU}_{1,n}^{ns}$ the open subset of
$\pi^*\OO(1)\ot\chi$-semistable points.

\begin{thm}\label{GIT-thm} 
(i) Let $(C,p_\bullet,v_\bullet)$ be in $\wt{\UU}_{1,n}^{ns}$, and let $E\sub C$ be the minimal elliptic subcurve.
Let $I\sub\{1,\ldots,n\}$ be the set of all $i$ such that $p_i\in E$, $J\sub \{1,\ldots,n\}$ the set of all $j$ such that 
$p_j\not\in E$ and there are at least $3$ distinguished points on the irreducible component of $C$ containing $p_j$.
Finally define $I_0\sub I$ by
$$I_0=\begin{cases} \emptyset, & E \text{ is either smooth or nodal },\\ \{i\in I \ |\ p_i\in E'\sub E, 
E' \text{ has only } 2 \text{ special points}\},
& \text{otherwise.}\end{cases}$$
Then for a rational character $\chi=\sum a_i \be_i$ of $\G_m$, the point $(C,p_\bullet,v_\bullet)$ is in $\wt{\UU}_{1,n}^{ns}(\chi)$
 if and only if 
\begin{itemize}
\item $a_i\ge 0$ for all $i$, and $a_i=0$ for $i\not\in I\cup J$;
\item $\sum_{i\in I} a_i\ge 1$;
\item $\sum_{i\in I_0} a_i\le 1$.
\end{itemize}
In particular, $\wt{\UU}_{1,n}^{ns}(\chi)$ is empty unless all $a_i\ge 0$ and $\sum_{i=1}^n a_i\ge 1$.

\smallskip

\noindent
(ii) With the notation of (i), $(C,p_\bullet,v_\bullet)$ has a finite reduced stabilizer subgroup in $\G_m^n$ 
if and only if $I\cup J=\{1,\ldots,n\}$ and $I_0\neq I$. If this is the case, and in addition $\chi=\sum a_i\be_i$ satisfies
$a_i>0$, $\sum_{i\in I}a_i>1$ and $\sum_{i\in I_0}a_i<1$,
then $(C,p_\bullet,v_\bullet)$ is $\pi^*\OO(1)\ot\chi$-stable. 

\smallskip

\noindent
(iii) Let us consider the following set of hyperplanes (walls) in $\R^n$: $a_i=0$ and $\sum_{i\in S} a_i=1$,
for all $i=1,\ldots,n$ and all nonempty subsets $S\sub [1,n]$. Let $\WW\sub\R^n$ be the union of these walls.
Then for $\chi=\sum a_i\be_i$ in a fixed connected component (chamber) 
$\Si$ of $\R^n\setminus\WW$ the stability coincides with
semistability and depends only on $\Si$.
\end{thm}

\Pf . (i) By Lemma \ref{cones-lem2}(ii,iii), we have
$$\Om_1(C,p_\bullet)=\{\be_i \ |\ i\in I\}, \ \ \bC(\Om_0(C,p_\bullet))=\bC(\be_i \ |\  i\in J\cup I\setminus I_0).$$
It is easy to see that for subsets of indices $T\sub S$ one has
$$\Conv(\be_i \ | i\in S)+\bC(\be_i \ |\ i\in S\setminus T)=\{ \sum_{i\in S} x_i\be_i \ |\ x_i\ge 0, \sum_{i\in S} x_i\ge 1, \sum_{i\in T} x_i\le 1\}.
$$
Thus, by Lemma \ref{cones-lem1}, we have
$$\Conv(C,p_\bullet)=\{ \sum_{i\in I\cup J} x_i\be_i \ |\ x_i\ge 0, \sum_{i\in I} x_i\ge 1, \sum_{i\in I_0} x_i\le 1\}.$$
Hence, we can rewrite the condition $\chi\in \Conv(C,p_\bullet)$ in the stated form.

\noindent
(ii) The first assertion follows immediately from Lemma \ref{Aut-curve-lem}. For the second we observe that if $\chi$
satisfies these strict inequalities then there exist positive rational numbers $(t_i)_{i\in I}$
and $(s_i)_{i\in J\cup I\setminus I_0}$, with $\sum_i t_i=1$, such that 
$$\chi=\sum_{i\in I} t_i\be_i+ \sum_{i\in J\cup I\setminus I_0} s_i\be_i.$$
Then for sufficiently divisible $N>0$ there exist functions $f_i\in A$, for $i\in J\cup I\setminus I_0$, such that
$$f=\prod_{i\in I}x_i^{Nt_i}\cdot \prod_{i\in J\cup I\setminus I_0} f_i^{Ns_i}$$
is an invariant section of $\OO(N)\ot\chi^N$ that does not vanish at $(C,p_\bullet,v_\bullet)$.
Then for any point $(C',p'_\bullet,v'_\bullet)$ of the open set $(f\neq 0)$ we have $I'\supset I$ and 
$J'\cup I'\setminus I'_0\supset J\cup I\setminus I_0$ (where $(I',I'_0,J')$ are defined for $(C',p'_\bullet)$ in the same
way as $(I,I_0,J)$ for $(C,p_\bullet)$). Hence, we get that every point of $(f\neq 0)$ has a finite stabilizer subgroup in 
$\G_m^n$, and hence $\pi^*\OO(1)\ot\chi$-stable.

\noindent
(iii) The fact that the stability condition depends only on a chamber is immediate from (i). By (ii), it remains to show that for
$\chi\not\in\WW$ all $\pi^*\OO(1)\ot\chi$-semistable points have 
$I\cup J=\{1,\ldots,n\}$ and $I_0\neq I$. But this follows again from (i): if $i\not\in I\cup J$
then $a_i=0$, and if $I_0=I$ then $\sum_{i\in I}a_i=1$.
\ed


In particular, we can describe for which characters we get modular compactifications of $\MM_{1,n}$ (over $\Z[1/6]$), i.e.,
(nonempty) 
open substacks in the stack $\UU_{1,n}$ of all pointed curves of arithmetic genus $1$, which are proper over $\Z[1/6]$.

\begin{cor}\label{mod-comp-cor}
For every $\chi=\sum a_i\be_i$ with $a_i>0$, $\sum_{i=1}^n a_i>1$ and $\chi\not\in\WW$, the quotient-stack
$$\UU_{1,n}^{ns}(\chi):=\wt{\UU}_{1,n}^{ns}(\chi)/\G_m^n$$
contains $\MM_{1,n}$ as an open substack and is a proper Deligne-Mumford stack over $\Z[1/6]$, 
with the projective coarse moduli space $\wt{\UU}^{ns}_{1,n}\sslash_\chi \G_m^n$. Hence, $\UU_{1,n}^{ns}(\chi)$
is a modular compactification of $\MM_{1,n}$ over $\Z[1/6]$.
\end{cor}

\Pf . Theorem \ref{GIT-thm}(ii) immediately implies that every $(C,p_\bullet,v_\bullet)$ with smooth $C$ is
$\pi^*\OO(1)\ot\chi$-stable. Hence, $\UU_{1,n}^{ns}(\chi)$ contains $\MM_{1,n}$ as an open substack.
Note also that by Theorem \ref{GIT-thm}(ii,iii), all points in $\wt{\UU}_{1,n}(\chi)$ are stable and have finite reduced
stabilizers. This implies that $\UU_{1,n}^{ns}(\chi)$ is a Deligne-Mumford stack and that
$\wt{\UU}^{ns}_{1,n}\sslash_\chi \G_m^n$ is its coarse moduli space (by \cite[Prop.\ 7.7]{Alper}).
Since $\wt{\UU}^{ns}_{1,n}\sslash_\chi \G_m^n$ is projective (by Proposition \ref{GIT-proj-prop}), it follows that
 $\UU_{1,n}(\chi)$ is proper (see the proof of \cite[Cor.\ 2.4.4]{P-ainf}).
\ed

The following particular case of Theorem \ref{GIT-thm}(i) 
will be useful in comparing the semistable loci in $\wt{\UU}_{1,n}^{ns}$ with
the Smyth's $m$-stable curves.

\begin{cor}\label{GIT-cor} 
Let $(C,p_\bullet,v_\bullet)\in \wt{\UU}_{1,n}^{ns}$, and let $\chi=\sum a_i\be_i$ with $a_i\ge 0$ 
and $\sum_{i=1}^n a_i\ge 1$.
Assume that $C$ coincides with its minimal elliptic subcurve $E$. Then the point $(C,p_\bullet,v_\bullet)$ is $\pi^*\OO(1)\ot\chi$-semistable
if and only if $\sum_{i\in I_0} a_i\le 1$.
\end{cor}


Recall that Smyth defines in \cite{Smyth-modular} a modular compactification $\ov{\MM}_{g,n}(\ZZ)$
of $\MM_{g,n}$ for each {\it extremal assignment} $\ZZ$. Roughly speaking, this is a combinatorial rule that assigns
to each stable curve some of its irreducible components, in a way compatible with specialization.
The moduli stack $\ov{\MM}_{g,n}(\ZZ)$ then classifies {\it $\ZZ$-stable curves} defined as contractions of the $\ZZ$-assigned components in usual stable curves (with some natural conditions on the obtained singular points, see 
\cite[Def.\ 1.8]{Smyth-modular}).

We are interested in one particular extremal assignment $\ZZ^u$ which assigns to each stable curve all of its unmarked components, see \cite[Ex.\ 1.12]{Smyth-modular}. It turns out that
in the chamber containing $\chi$ with large positive $a_i$ our GIT stability condition gives precisely
the $\ZZ^u$-stability.

\begin{prop}\label{new-stability-prop} 
Let $\chi=\sum a_i\be_i$ be any rational character with $a_i>1$ for all $i$.
Then for $(C,p_\bullet,v_\bullet)\in \wt{\UU}_{1,n}^{ns}$ the following conditions are equivalent:
\begin{enumerate}
\item $(C,p_\bullet,v_\bullet)$ is $\pi^*\OO(1)\ot\chi$-semistable;
\item $(C,p_\bullet,v_\bullet)$ is $\pi^*\OO(1)\ot\chi$-stable;
\item there are $\ge 3$ special points on the normalization of every rational irreducible component of $C$;
\item $(C,p_\bullet)$ is $\ZZ^u$-stable.
\end{enumerate}
\end{prop}

\Pf . As before, we denote by $E$ the minimal elliptic subcurve of $C$.
By Theorem \ref{GIT-thm}, in our situation the semistability is equivalent to stability and corresponds to the
conditions $I\cup J=[1,n]$ and $I_0=\emptyset$. By definition, this means that every rational irreducible component of $C$
has at least $3$ special points, and in addition, if $E$ is the elliptic $m$-fold curve then every irreducible component of $E$
also has at least $3$ special points, which shows the equivalence of (1) and (2) with (3). 
Finally, to prove the equivalence with (4), we observe that if $(C,p_\bullet)$ is $\ZZ^u$-stable then
$C$ is obtained from a stable curve $\wt{C}$ by contracting the components without marked points. Then every rational component $R$ of $C$ is the image of a rational component
$\wt{R}$ of $\wt{C}$ such that the map $\wt{R}\to R$ is a bijection. Thus, the images of $\ge 3$ special points on $\wt{R}$
give $\ge 3$ special points on $R$. Conversely, assume $C$ satisfies condition (3). If $C$ is at most nodal then $(C,p_\bullet)$ is
stable, so it is also $\ZZ^u$-stable (since there is at least one marked point on each irreducible component by the definition
of $\wt{\UU}_{1,n}^{ns}$). Now assume that $E$ is the elliptic $m$-fold curve. Then $E=\cup_{i=1}^m E_i$ is obtained 
as the contraction of the nodal curve $\wt{E}=\cup_{i=0}^m E_i$ with $m+1$ components, where $E_0$ is an elliptic curve
with no marked points, $E_1, \ldots, E_m$ are attached to $E_0$ transversally at $m$ distinct points. Now $C$ is the union
of $E$ and some rational tails (with $l$-fold rational singularities), and we can define $\wt{C}$ to be the union of $\wt{E}$ and
these rational tails. More precisely, if $R$ is a rational tail attached transversally to a smooth point of $E$,
which lies in $E_i$, then we attach $R$ to $E_i\sub\wt{E}$ in the same way. On the other hand, if $R$ is attached to the singular
point of $E$ then we attach $R$ to $E_0$ transversally at a point distinct from all the other special points on $E_0$.
Then $\wt{C}$ with the induced marked points will be a stable curve, proving that $(C,p_\bullet)$ is $\ZZ^u$-stable.
\ed

\begin{cor}\label{coarse-Z-cor} 
The coarse moduli space $\ov{M}_{1,n}(\ZZ^u)$ is projective over $\Spec\Z[1/6]$. 
\end{cor}

\begin{rems} 1. Note that all the modular compactifications of $\MM_{1,n}$
considered in Corollary \ref{mod-comp-cor}
are {\it semistable} in the sense of \cite[Def.\ 1.2]{Smyth-modular}, i.e., the normalization of every
rational component contains at least two distinguished points 
(in our case every rational component contains a marked point and at least one singular point of $C$).
The only case when $\UU_{1,n}^{ns}(\chi)$ is a {\it stable} modular compactification (i.e., the normalization of every
rational component contains at least three distinguished points) is the one considered in
Proposition \ref{new-stability-prop}.

\noindent
2. Let us work over $\C$.
Recall that {\it Kontsevich's compactification} of $M_{g,n}$ is the topological space $K\ov{M}_{g,n}$ obtained
as a certain quotient of $\ov{M}_{g,n}$ (see \cite{Kon}, \cite{Looijenga}). Namely, two stable pointed curves $C_1$ and $C_2$
are identified if the corresponding curves $\ov{C}_1$ and $\ov{C}_2$, obtained by contracting all unmarked components,
are isomorphic, and this isomorphism preserves the extra data about the arithmetic genus of each 
contracted connected subcurve in the original curves. In the case $g=0$ the spaces $K\ov{M}_{0,n}$ have been realized algebraically by Boggi \cite{Boggi}, and these spaces coincide with the GIT quotients $\wt{\UU}_{0,n}\sslash_\chi \G_m^n$
for appropriate characters $\chi$ (see \cite[Sec.\ 5]{P-ainf}). The moduli spaces $\ov{M}_{1,n}(\ZZ^u)$,
realized in Proposition \ref{new-stability-prop} as GIT quotients, are closely related to $K\ov{M}_{1,n}$, at least for small $n$.
Indeed, we have a natural birational map 
\begin{equation}\label{M-Z-map}
\ov{M}_{1,n}\to \ov{M}_{1,n}(\ZZ^u).
\end{equation}
It is easy to see that 
in the cases $n=2$ and $n=3$ this map is regular and the fibers are precisely the equivalence classes defining
$K\ov{M}_{1,n}$, so in these cases $\ov{M}_{1,n}(\ZZ^u)$ is an algebraic realization of $K\ov{M}_{1,n}$
(more precisely, \eqref{M-Z-map} blows down the divisor of curves with the unmarked component of genus $1$ to the locus
in $\ov{M}_{1,n}(\ZZ^u)$ corresponding to curves with a cusp). In the case $n=4$
the rational morphism \eqref{M-Z-map} is not defined at the elliptic bridge $C=E\cup R_1\cup R_2$,
where $E$ is the unmarked elliptic curve, the rational components $R_1$ and $R_2$ are attached to $E$ and
both are equipped with two marked points $p_1,p_2\in R_1$, $p_3,p_4\in R_2$. 
However, using \cite[Lem.\ 4.2]{Smyth-II}, we can fix this by identifying
the following three points of $\ov{M}_{1,4}(\ZZ^u)$: 1) the curve $E\cup R$, where $E$ is cuspidal, $R$ is attached to the
singular point of $E$, $p_1,p_2\in E$, $p_3,p_4\in R$;
2) the same curve but with $p_1,p_2\in R$, $p_3,p_4\in E$; 3) the tacnode $E$ with $p_1,p_2$ on one component and 
$p_3,p_4$ on another component. If we similarly identify points in two other triples of points in $\ov{M}_{1,4}(\ZZ^u)$, corresponding to different partitioning of the marked points into two pairs, the map \eqref{M-Z-map} will induce a regular morphism to the resulting
quotient, which will be homeomorphic to $K\ov{M}_{1,4}$ (and thus, will give an algebraic model of it, since identifying of a finite number of points can be done algebraically, using an affine open containing the finite set of points to be identified).
\end{rems}



\subsection{Connection with Smyth's $m$-stable curves}\label{m-stable-sec}

We continue to work over $\Z[1/6]$.
Recall that $\ov{\MM}_{1,n}(m)$ denotes the moduli stack of Smyth's $m$-stable curves, and
$\ov{M}_{1,n}(m)$ is its coarse moduli stack (see \cite{Smyth-I}).
Elsewhere we showed that for each $m$, $(n-1)/2\le m\le n-1$, 
there is a natural morphism
$$\rho:\ov{\MM}_{1,n}(m)\to \UU_{1,n}^{ns},$$
associating with $(C,p_\bullet)$ the curve $(\ov{C},p_\bullet)$, where $\ov{C}$ is obtained from $C$ by contracting
the components without the marked points (see \cite[Prop.\ 1.5.1]{P-ainf-more-pts}).

Actually, in the case $m=n-1$ every $(n-1)$-stable curve is already in $\UU_{1,n}^{ns}$, so in this case
the morphism $\rho$ is the inclusion of an open substack 
(see \cite[Prop.\ 3.3.1]{P-krich} and \cite[Thm.\ 1.5.7]{LP}).
However, for all $m<n-1$ there exist $m$-stable curves with unmarked components.
We are going to show in Proposition \ref{Smyth-GIT-prop} below
that for $m=n-1$, $n-2$ and $n-3$, the image of $\rho$ is contained in
some $\pi^*\OO(1)\ot\chi$-semistable loci.

\begin{lem}\label{orbit-closure-lem} 
Assume that $(C',p'_\bullet,v'_\bullet)\in \wt{\UU}_{1,n}^{ns}$ is in the closure of the
$\G_m^n$-orbit of $(C,p_\bullet,v_\bullet)\in \wt{\UU}_{1,n}^{ns}$.

\noindent
(i) Assume that for some $i\neq j$ the points $p'_i$ and $p'_j$ belong to the same component of a rational
tail of $C'$. Then 
$p_i$ and $p_j$ belong to the same component of $C$.

\noindent
(ii) Assume now that $p'_i$ and $p'_j$ belong to the same rational tail of $C'$, and that $C$ coincides with its minimal
elliptic subcurve.
Then $p_i$ and $p_j$ belong to the same component of $C$.
\end{lem}

\Pf . (i) Pick $m$ such that $p'_m$ belongs to the minimal elliptic subcurve in $C'$. Then for $(C',p'_\bullet,v'_\bullet)$
we have $x_m\neq 0$ (by Corollary \ref{H1-cor}) and 
$c_{ij}(m)=h_{mi}(p'_j)\neq 0$, $c_{ji}(m)\neq 0$ (see Lemma \ref{f-h-h-lem}(iii)). 
Hence, the same inequalities should hold for $(C,p_\bullet,v_\bullet)$. In particular, $p_m$ belongs to the minimal
elliptic subcurve $E\sub C$.
Assume that $p_i$ and $p_j$ lie on different components
of $C$. If one of them lies on a rational tail, then from Lemma \ref{f-h-h-lem}(iii) we get that
either $c_{ij}(m)=0$ or $c_{ji}(m)=0$. Hence, both $p_i$ and $p_j$ lie on $E$. Then by Lemma \ref{f-h-h-lem}(i),
the nonvanishing of $c_{ij}(m)$ implies that one of the points $(p_i,p_j)$, say $p_i$, lies on the same component as $p_m$.
Now let us consider
the $\G_m^n$-invariant rational function $a_{mi}^2b_{mi}/c_{ij}(m)^2$ on $\wt{\UU}_{1,n}(m)$.
It is well defined on $(C,p_\bullet,v_\bullet)$ and on $(C',p'_\bullet,v'_\bullet)$, so it has to take the same values. But it vanishes at 
$(C',p'_\bullet,v'_\bullet)$ (since $a_{mi}$ vanishes), and one can check using Lemma \ref{f-h-h-lem}(i) that it does not
vanish on $(C,p_\bullet,v_\bullet)$. This contradiction shows that $p_i$ and $p_j$ belong to the same component of $C$.

\noindent
(ii) Pick $m$ as in (i). Then for $(C',p'_\bullet,v'_\bullet)$ we have $x_m\neq 0$ and either $c_{ij}(m)\neq 0$ or $c_{ji}(m)\neq 0$.
We can assume that $c_{ij}(m)\neq 0$. Then we also have the same inequality for $(C,p_\bullet,v_\bullet)$.
Assume that $p_i$ and $p_j$ lie on different components of $C=E$.
By Lemma \ref{f-h-h-lem}(i), this implies that either $p_i$ belongs to the same component as $p_m$, or
$p_j$ belongs to the same component as $p_m$. In the former case the $\G_m^n$-invariant function
$a_{mi}^2b_{mi}/c_{ij}(m)^2$ does not vanish at $(C,p_\bullet,v_\bullet)$, so we can finish the proof as in (i).
If $p_j$ and $p_m$ lie on the same component then 
the same argument works with the $\G_m^n$-invariant function $a_{mi}^2b_{mj}/c_{ij}(m)^2$.
\ed

\begin{prop}\label{Smyth-GIT-prop} Consider a rational character $\chi=\sum_{i=1}^n a_i\be_i$
of $\G_m^n$, where $a_i\ge 0$ for all $i$, and let $\UU_{1,n}^{ns}(\chi)\sub \UU_{1,n}^{ns}$ denote the quotient-stack
of the $\pi^*\OO(1)\ot\chi$-semistable locus $\wt{\UU}_{1,n}^{ns}(\chi)$ by $\G_m^n$.

\noindent
(i) For $n\ge 2$ assume that 
\begin{equation}\label{n-1-stability-eq}
\sum_{i=1}^n a_i\ge 1, \ \text{ and }\ \  \sum_{i\in S} a_i\le 1 \ \ \forall S\sub [1,n], |S|=n-2.
\end{equation}
Then we have an inclusion of an open substack
$$\ov{\MM}_{1,n}(n-1)\simeq\rho(\ov{\MM}_{1,n}(n-1))\sub \UU_{1,n}^{ns}(\chi).$$
Moreover, if $a_i>0$ for all $i$ and all the inequalities in \eqref{n-1-stability-eq} are replaced by the strict ones
then the above inclusion becomes an equality, and all points of $\wt{\UU}_{1,n}^{ns}(\chi)$ are
$\pi^*\OO(1)\ot\chi$-stable. In this case $\rho$ induces an isomorphism
$$\ov{M}_{1,n}(n-1)\simeq \wt{\UU}_{1,n}^{ns}\sslash_\chi \G_m^n$$

\noindent
(i') Assume that $\chi$ satisfies \eqref{n-1-stability-eq} and in addition, $\sum_{i=1}^n a_i>1$ and $a_i>0$ for all $i$.
Then $\rho$ induces an isomorphism
$$\ov{M}_{1,n}(n-1)\simeq (\wt{\UU}_{1,n}^{ns}\sslash_\chi \G_m^n)^{\nu},$$
where for an irreducible scheme $X$ we denote by $X^\nu$ the normalization of $X_{\red}$, and the morphism 
$\ov{M}_{1,n}(n-1)\to \wt{\UU}_{1,n}^{ns}\sslash_\chi \G_m^n$ is a bijection.

\noindent
(ii) For $n\ge 3$ assume that 
\begin{equation}\label{n-2-stability-eq}
\sum_{i\in S} a_i\ge 1 \ \ \forall S\sub [1,n], |S|=n-2,  \ \text{ and }\ \  \sum_{i\in S'} a_i\le 1 \ \ \forall S'\sub [1,n], |S'|=n-3.
\end{equation}
Then we have an inclusion
$$\rho(\ov{\MM}_{1,n}(n-2))\sub \UU_{1,n}^{ns}(\chi).$$
Moreover, if $a_i>0$ for all $i$ and all the inequalities in \eqref{n-1-stability-eq} are replaced by the strict ones
then the above inclusion becomes an equality, and all points of $\UU_{1,n}^{ns}(\chi)$ are
$\pi^*\OO(1)\ot\chi$-stable. 

\noindent
(ii') Assume that $\chi$ satisfies \eqref{n-2-stability-eq} and in addition, $a_i>0$ for all $i$ and
$\sum_{i\in S} a_i>1$ for all $S$ with $|S|=n-2$. 
Then $\rho$ induces an isomorphism
$$\ov{M}_{1,n}(n-2)\simeq (\wt{\UU}_{1,n}^{ns}\sslash_\chi \G_m^n)^{\nu},$$
and the morphism 
$\ov{M}_{1,n}(n-2)\to \wt{\UU}_{1,n}^{ns}\sslash_\chi \G_m^n$ is a bijection.

\noindent
(iii) For $n>4$ assume that $a_i=\frac{1}{n-4}$ for all $i$. Then we have an inclusion
$$\rho(\ov{\MM}_{1,n}(n-3))\sub \UU_{1,n}^{ns}(\chi).$$
\end{prop}

\Pf . (i) Let $(C,p_\bullet)$ be an $(n-1)$-stable curve, i.e., we have $C=E$, $(C,p_\bullet)$ has no infinitesimal symmetries, and $E$ is either at most nodal, or the elliptic $m$-fold curve with $m\le n-1$.

We claim that $|I_0|\le n-2$. Indeed, this is clear 
if $C$ is at most nodal. On the other hand, if $C=E$ is the elliptic $m$-fold curve with $m\le n-1$ then $C$ has
at most $m-1\le n-2$ components with a unique marked point, i.e., $|I_0|\le n-2$, as claimed.
Thus, by Corollary \ref{GIT-cor}, the $\pi^*\OO(1)\ot\chi$-semistability of $(C,p_\bullet)$ follows from 
$\sum_{i\in I_0} a_i\le 1$.

Next, assume that $a_i>0$ and the strict inequalities in \eqref{n-1-stability-eq} hold, and let $(C,p_\bullet,v_\bullet)$
be $\pi^*\OO(1)\ot\chi$-semistable. Then from the criterion
of Theorem \ref{GIT-thm}(i) we obtain that $|I|>n-2$. Also, by Theorem \ref{GIT-thm}(ii), the point
$(C,p_\bullet,v_\bullet)$ is actually $\pi^*\OO(1)\ot\chi$-stable and $(C,p_\bullet)$ has no nontrivial infinitesimal symmetries.
Hence, any rational tail of $C$ should contain $\ge 2$ marked points, and we get $|I|=n$, i.e., $C=E$. 
Note also that if $C$ is the elliptic $m$-fold curve then we necessarily have $m\le n-1$: if $m=n$ then
$C$ would have one marked point on every component, so $(C,p_\bullet)$ would have $\G_m$ as the group of automorphisms.

Thus, in this case we have $\ov{\MM}_{1,n}(n-1)\simeq \UU_{1,n}^{ns}(\chi)$.
Furthermore, by Corollary \ref{mod-comp-cor}, the corresponding GIT quotient
is the coarse moduli space of the stack $\UU_{1,n}^{ns}(\chi)$.
Hence, it is isomorphic to $\ov{M}_{1,n}(n-1)$.

\noindent
(i') First, we observe that by part (i), we have an inclusion
$$\ov{\MM}_{1,n}(n-1)\sub \UU_{1,n}^{ns}(\chi).$$

Next, using \eqref{n-1-stability-eq}, the positivity of $a_i$'s and Theorem \ref{GIT-thm}(i), we see that any 
$(C,p_\bullet,v_\bullet)$ in $\wt{\UU}_{1,n}^{ns}(\chi)$ has $I\cup J=[1,n]$ and $|I|\ge n-2$. Hence, either
$C=E$ or $C=E\cup R$, where $R\simeq\P^1$, and $R$ contains two marked points.
In the former case the condition $\sum_{i=1}^n a_i>1$ gives that $|I_0|<n$, so $E$ cannot be the elliptic $n$-fold
curve. Thus, by Theorem \ref{GIT-thm}(ii), the only case when $(C,p_\bullet,v_\bullet)$ has an infinite
stabilizer is when $|I_0|=|I|=n-2$, so $n\ge 3$, $C=E\cup R$, where $E$ is the elliptic $(n-2)$-fold curve with one marked
point on each component, and $R$ is attached to the singular point of $E$. Note that we have such pointed curves 
$C_{ij}\in \UU_{1,n}^{ns}(\chi)$, for each pair of indices $i<j$, corresponding to the marked points on $R$.
Note also that the corresponding $\G_m^n$-orbits in $\wt{\UU}_{1,n}^{ns}(\chi)$ are $(n-1)$-dimensional. 

By Theorem \ref{GIT-thm}(ii), 
all the points in $\ov{\MM}_{1,n}(n-1)$ with $|I_0|<n-2$ are $\pi^*\OO(1)\ot\chi$-stable, so they correspond to closed orbits in
$\wt{\UU}_{1,n}^{ns}(\chi)$. Assume that $(C,p_\bullet)\in \ov{\MM}_{1,n}(n-1)$ has $|I_0|\ge n-2$.
Then $|I_0|=n-2$ and $C=E$ is the elliptic $(n-1)$-fold curve, with two marked points
on one component and one marked point on each of the remaining $(n-2)$ components.
Furthermore, if $C_{ij}$ is in the closure of the $\G_m^n$-orbit corresponding to
such a curve $(C,p_\bullet)\in \ov{\MM}_{1,n}(n-1)$ then
by Lemma \ref{orbit-closure-lem}(i), we see that $p_i$ and $p_j$ lie on one component of $C$.
 

Thus, there is at most one curve $(C,p_\bullet)\in \ov{\MM}_{1,n}(n-1)$ that has $C_{ij}$ in the closure of the corresponding orbit. 
It follows that the natural birational morphism
$$\ov{M}_{1,n}(n-1)\to \wt{\UU}_{1,n}^{ns}\sslash_\chi \G_m^n$$
is injective on geometric points. Since $\ov{M}_{1,n}(n-1)$ is normal and projective over $\Z[1/6]$ (see \cite[Cor.\ 1.5.15]{LP}),
this implies that it is a bijection and coincides with the normalization of $\wt{\UU}_{1,n}^{ns}\sslash_\chi \G_m^n$.

\noindent
(ii) Assume that $\chi$ satisfies \eqref{n-2-stability-eq}.
If $(C,p_\bullet)$ is an $(n-2)$-stable curve then either $C=E$, or
$C=E\cup R$ for one rational component $R$, containing $2$ marked points. Let $C\to \ov{C}$ be
the contraction of the unmarked components. Note that if $C=E$ then every component of $C$ has at least one marked
point (this follows from the condition that $C$ has no infinitesimal symmetries), so we have $\ov{C}=C$. 
In the case when $C=E\cup R$,
the minimal elliptic subcurve $E$ can contain one unmarked component, to which $R$ is attached, so in this case 
$\ov{C}=\ov{E}\cup R$, where $\ov{E}$ is the image of $E$.
Thus, in any case for the curve $(\ov{C},p_\bullet)$ we have $|I|\ge n-2$. 
In the case when $E$ is at most nodal, $\ov{E}$ is also at most nodal, so we deduce that it is $\pi^*\OO(1)\ot\chi$-semistable. 
If $E$ is the elliptic $m$-fold curve with $m\le n-2$ then either $\ov{E}=E$ or $\ov{E}$ is the elliptic $(m-1)$-fold curve,
so we have $|I_0|\le n-3$. Thus, in the case $C=E$
the $\pi^*\OO(1)\ot\chi$-semistability follows from Corollary \ref{GIT-cor}. In the case $C=E\cup R$, we have to
observe in addition that $I\cup J=\{1,\ldots,n\}$, since $R$ has three special points.

Now assume that $a_i>0$ and the strict inequalities in \eqref{n-2-stability-eq} hold, and let $(C,p_\bullet,v_\bullet)$
be a $\pi^*\OO(1)\ot\chi$-semistable point. Then by Theorem \ref{GIT-thm}(ii), this point is stable and $(C,p_\bullet)$
has no infinitesimal symmetries. From the criterion of Theorem \ref{GIT-thm}(i) we derive that $|I|>n-3$. Hence,
either $|I|=n$, or $|I|=n-2$ and $C=E\cup R$, where $R$ is a smooth rational component with two marked points on it.

In the case when either $C=E$ or $R$ is attached to a smooth point on $E$ we claim 
that $(C,p_\bullet)$ is itself $(n-2)$-stable. For this we have to show that the singularities of $C$
are either nodes or elliptic $m$-fold points
with $m\le n-2$. Thus, in the case $C=E$ we have to rule out the case when $E$ is the elliptic $m$-fold curve
with $m\ge n-1$. But in this case the condition that $E$ has no infinitesimal symmetries implies that $m=n-1$ and
exactly one component of $E$ contains two marked points, and the remaining components
contain one marked point each. Hence, we would have $|I_0|=m-1=n-2$, and the condition
$\sum_{i\in I_0} a_i\le 1$ would contradict the strict version of \eqref{n-2-stability-eq}.
Now assume $C=E\cup R$, where $E$ is the elliptic
$m$-fold curve, $R\simeq \P^1$ contains two marked points (so $E$ contains $n-2$ marked points), and the intersection
of $E$ and $R$ is nodal. Since $C$ has no infinitesimal symmetries we derive that $m\le n-2$, which implies
that $(C,p_\bullet)$ is $(n-2)$-stable. 

It remains to consider the case when $C=E\cup R$, where $R\simeq \P^1$ contains two marked points, and $R$
is attached to a singular point $q\in E$. Assume first that $E$ is the standard $l$-gon. Then 
let $\wt{C}=\wt{E}\cup R$, where $\wt{E}$ is the standard
$(l+1)$-gon obtained by gluing in the extra component instead of $q$ in $E$, and $R$ is attached to a smooth point
on this extra component. Then $(C,p_\bullet)=\rho(\wt{C},p_\bullet)$, and $(\wt{C},p_\bullet)$ is $(n-2)$-stable.
Now assume that $E$ is the elliptic $m$-fold curve. Then the condition that $C$ has no infinitesimal symmetries implies that 
$m\le n-3$. Thus, we have $(C,p_\bullet)=\rho(\wt{C},p_\bullet)$, where $\wt{C}=\wt{E}\cup R$, with $\wt{E}$ being
the elliptic $(m+1)$-fold curve, with one unmarked component to which $R$ is attached. Since $m+1\le n-2$,
the curve $(\wt{C},p_\bullet)$ is $(n-2)$-stable.

\noindent
(ii') Our work in (ii) shows that the morphism 
$$\ov{\MM}_{1,n}(n-2)\to \UU_{1,n}^{ns}(\chi),$$
induced by $\rho$, is an injection on geometric points. Let us consider the induced birational morphism 
\begin{equation}\label{coarse-(n-2)-morphism}
\ov{M}_{1,n}(n-2)\to \wt{\UU}_{1,n}^{ns}\sslash_\chi \G_m^n.
\end{equation}
Since $\ov{M}_{1,n}(n-2)$ is normal and projective, as in (i'), it is enough to check that this morphism is injective
on geometric points.

Using Theorem \ref{GIT-thm}(i) we see that our assumptions on $\chi$ imply that
every $(C,p_\bullet,p_\bullet)\in\wt{\UU}_{1,n}^{ns}(\chi)$ has $I\cup J=[1,n]$, $|I|\ge n-3$ and
$|I_0|\le n-3$. Hence, either $C=E$, or $C=E\cup R$, where $R$ is a rational tail with $\le 3$ marked points.
In the former case the condition $|I_0|\le n-3$ imples that $E$ is either at most nodal, or the elliptic $m$-fold curve
with $m\le n-2$. The only case when $(C,p_\bullet,v_\bullet)$ has an infinite stabilizer is when
$|I_0|=|I|=n-3$ (in particular, $n\ge 4$), 
so $C=E\cup R$, where $E$ is the elliptic $(n-3)$-fold curve with one marked point on each component,
and $R$ is attached to the singular point $q$ of $E$. There are two types of such pointed curves: $C_{ijk}(t)$, for $i<j<k$,
corresponds to the case when $R=\P^1$ and $p_i,p_j,p_k\in R$, so that the points $(q,p_i,p_j,p_k)$ have the cross-ratio $t$;
and $C_{i;jk}$, where $R$ is the union of two components, $p_i,p_j,p_k\in R$, and $p_i$ belongs to the component
attached to $E$. In any case the corresponding $\G_m^n$-orbit in $\wt{\UU}_{1,n}^{ns}(\chi)$ is $(n-1)$-dimensional.
 
Now we have to determine which $\G_m^n$-orbits in $\wt{\UU}_{1,n}^{ns}(\chi)$ contain the above pointed curves in their closure.
Recall that we can assume that $n\ge 4$ (otherwise, all orbits are closed).
Note that by Theorem \ref{GIT-thm}(ii),
the curves in $\rho(\ov{\MM}_{1,n}(n-2))$ with $|I_0|<n-3$ are stable (since they have $|I|\ge n-2$), 
so the corresponding $\G_m^n$-orbits
are closed in $\wt{\UU}_{1,n}^{ns}(\chi)$. Otherwise, we have $|I_0|=n-3$, so $|I|$ is either $n-2$ or $n$.
There are two types of such curves in $\rho(\ov{\MM}_{1,n}(n-2))$: 
\begin{enumerate}
\item
$C=E\cup R$, where $E$ is the elliptic $(n-2)$-fold curve, with one marked point on each component, and 
$R=\P^1$ (with two marked points) is attached to a smooth point of $E$;
\item
$C=E$, the elliptic $(n-2)$-fold curve, with $3$ marked points on one component and one on each of the remaining components. 
\end{enumerate}
 
Assume that we have $(C,p_\bullet)\in \rho(\ov{\MM}_{1,n}(n-2))$ that contains $C_{ijk}(t)$ in its closure.
Then $(C,p_\bullet)$ should be as in cases (1) or (2) above. In particular, $E\sub C$ is the elliptic $(n-2)$-fold curve.
Furthermore, by Lemma \ref{orbit-closure-lem}(i), $p_i$, $p_j$ and $p_k$ all lie on the same component of $C$,
hence, $(C,p_\bullet)$ is as in case (2). We claim that $t$ is equal to the cross-ratio
of $(q,p_i,p_j,p_k)$ on this component. Indeed, let $E_1\sub E=C$ be the component containing $p_i,p_j,p_k$, and 
let $p_m$ be a marked point on a different component of $C$ (recall that $n\ge 4$).
By Lemma \ref{f-h-h-lem}(i), the function $h_{km}$ is
a coordinate on $E_1$, equal to $0$ at $q$ and to $\infty$ at $p_k$. Hence, the cross-ratio of $(q,p_i,p_j,p_k)$ is given by
$$\frac{h_{km}(p_j)}{h_{km}(p_j)-h_{km}(p_i)}=\frac{c_{mj}(k)}{c_{mj}(k)-c_{mi}(k)}.$$
Now we observe that the latter ratio is $\G_m^n$-invariant, so it is preserved under passing to the closure of the 
$\G_m^n$-orbit. But for the limit curve $C_{ijk}(t)$ 
the restriction of $h_{km}$ to $R$ is still a coordinate on $R$, equal to $0$ at $q$
and to $\infty$ at $p_k$, so the limiting points on $R$ have the same cross-ratio.

Now assume that we have $(C,p_\bullet)\in \rho(\ov{\MM}_{1,n}(n-2))$ that contains $C_{i;jk}$ in its closure.
Then by Lemma \ref{orbit-closure-lem}(i), we get that $p_j$ and $p_k$ lie on the same component of $C$.
If we are in case (2) then applying Lemma \ref{orbit-closure-lem}(ii),
we get that $p_i$ belongs to the same component. But then
the ratio $c_{ij}(m)/c_{ik}(m)$ is preserved under passing to the closure, so we get a contradiction, since it is equal
to $1$ on $C_{i;jk}$, and it is not equal to $1$ on $C$. Hence, $C$ is as in case (1), with $p_j,p_k\in R$.
Furthermore, the condition $c_{ij}(m)\neq 0$ implies that $p_i$ belongs to the component to which $R$ is attached.

Thus, we showed that the morphism \eqref{coarse-(n-2)-morphism} is injective on geometric points as claimed.

\noindent
(iii) Let $(C,p_\bullet)$ be an $(n-3)$-stable curve. Then either $C=E$, or $C=E\cup R_1$, where $R_1$ contains $2$ or $3$
marked points (so $R_1$ can have $\le 2$ components), or $C=E\cup R_1\cup R_2$, where each $R_i$ contains $2$
marked points. Let us show that $\rho(C,p_\bullet)=(\ov{C},p_\bullet)$ is $\pi^*\OO(1)\ot\chi$-semistable.
Note that in all cases the curve $(\ov{C},p_\bullet)$ has $|I|\ge n-4$, so $\sum_{i\in I}a_i\ge 1$.
Note also that $(n-3)$-stable curves have no infinitesimal automorphisms, so we have $I\cup J=\{1,\ldots,n\}$.
Thus, we only have to check that $\sum_{I_0}a_i\le 1$, or equivalently, $|I_0|\le n-4$. 
Note the minimal elliptic subcurve $\ov{E}\sub \ov{C}$ is the image of $E\sub C$.
If $E$ is at most nodal then so is $\ov{E}$, hence $I_0=\emptyset$.
Now assume that $E$ is the elliptic $m$-fold curve. By definition, $(n-3)$-stability of $(C,p_\bullet)$
implies then that $m\le n-3$. Since $C$ has no infinitesimal automorphisms, at most $m-1$ components of $E$
contain a unique special point in addition to the singular point of $E$. Hence,
we deduce that $|I_0|\le m-1\le n-4$, as required.
\ed

\begin{rems}\label{Smyth-GIT-rem} 
1. It is not hard to see by the same methods as in Proposition \ref{Smyth-GIT-prop} that for $m<n-3$ there
is no $\chi$ for which all Smyth $m$-stable curves map to $\pi^*\OO(1)\ot\chi$-semistable points under $\rho$.
Let us illustrate this for $\chi=a\sum_i \be_i$. Let $m=n-k$.
Consider first a (smooth) elliptic curve $E$, and let $C=E\cup R_1\cup\ldots\cup R_{k-1}$, 
where each $R_i$ is a smooth rational component containing exactly two marked points, while
the remaining $n-2k+2$ marked points are on $E$.
Then this curve is both $m$-stable and belongs to $\UU_{1,n}^{ns}$, so it maps to itself by $\rho$. 
It also has $|I|=n-2k+2$, so we should have $a\ge 1/(n-2k+2)$ in order for it to be $\pi^*\OO(1)\ot\chi$-semistable.
On the other hand, let $C=E$ be the elliptic $m$-fold curve, and assume that $m-1$ components of $C$ 
each contain one marked point, while the last component contains the remaining $k+1$ marked points.
Then this curve is also $m$-stable and belongs to $\UU_{1,n}^{ns}$, and also has
$|I_0|=m-1=n-k-1$, so we should have $a\le 1/(n-k-1)$ in order for it to be  $\pi^*\OO(1)\ot\chi$-semistable. But 
this contradicts the previous inequality if $k>3$.

\noindent
2. Using Proposition \ref{Smyth-GIT-prop} we can explain the existence of a regular morphism 
\begin{equation}\label{(n-2)-(n-1)-morphism}
\ov{M}_{1,n}(n-2)\to \ov{M}_{1,n}(n-1)
\end{equation}
(see \cite[Cor.\ 4.5]{Smyth-II}). Namely,
let us take $\chi=\frac{1}{n-2}\cdot\sum \be_i$. Then by Proposition \ref{Smyth-GIT-prop}(ii), $\rho$ induces
a regular morphism
$$\ov{\MM}_{1,n}(n-2)\to \UU_{1,n}^{ns}(\chi)\to \wt{\UU}_{1,n}^{ns}\sslash_\chi \G_m^n,$$
which induces a morphism
$$\ov{M}_{1,n}(n-2)\to \wt{\UU}_{1,n}^{ns}\sslash_\chi \G_m^n.$$
Furthermore, since $\ov{M}_{1,n}(n-2)$ is normal (see \cite[Cor.\ 1.5.15]{LP}), it factors through a morphism
$$\ov{M}_{1,n}(n-2)\to (\wt{\UU}_{1,n}^{ns}\sslash_\chi \G_m^n)^{\nu}\simeq \ov{M}_{1,n}(n-1),$$
where we used Proposition \ref{Smyth-GIT-prop}(i').

\noindent
3. Using our work in the proof of Proposition \ref{Smyth-GIT-prop}, we can also easily see 
that the morphism \eqref{(n-2)-(n-1)-morphism} contracts each divisor $\De_{0,J}$ with $|J|=2$ to a point.
Here $\De_{0,J}$ consists of curves $E\cup R$, where $R\simeq\P^1$ contains two marked points in $J$.
By Lemma \ref{one-tail-orbit-lem}, the closure of the corresponding $\G_m^n$-orbit in $\wt{\UU}_{1,n}^{ns}$ contains
the curve $C_{ij}$, where $J=\{i,j\}$, from the proof of Proposition \ref{Smyth-GIT-prop}(i'). 
Hence, $\De_{0,J}$ gets contracted to the image of this point in 
$(\wt{\UU}_{1,n}^{ns}\sslash_\chi \G_m^n)^{\nu}\simeq \ov{M}_{1,n}(n-1)$,
where we take $\chi=\frac{1}{n-2}\cdot\sum \be_i$. Note that the corresponding points in $\ov{M}_{1,n}(n-1)$
are exactly those with $C$ being the elliptic $(n-1)$-fold curve.
Since each $\De_{0,J}$ is a Cartier divisor in $\ov{M}_{1,n}(n-1)$, it seems likely
that the morphism \eqref{(n-2)-(n-1)-morphism} identifies $\ov{M}_{1,n}(n-2)$ with the blow-up of
$\ov{M}_{1,n}(n-1)$ at the above ${n\choose 2}$ points (to prove this one has to compute the scheme-theoretic
preimage of each of these points under \eqref{(n-2)-(n-1)-morphism}).

\noindent
4. For $\chi=\frac{1}{n-4}\sum_{i=1}^n\be_i$ (i.e., in the situation of Proposition \ref{Smyth-GIT-prop}(iii)),
the map
$$\ov{M}_{1,n}(n-3)\to \wt{\UU}_{1,n}^{ns}\sslash_\chi \G_m^n,$$
induced by $\rho$, is still birational but it is not finite.
Namely, consider an $(n-3)$-stable curve of the form $C=E\cup R_1\cup R_2$, where $E$ is a standard $m$-gon with $m\ge 2$
and the tails $R_1$ and $R_2$ are attached to points $r_1\neq r_2$ of the same unmarked component $E_0\sub E$. 
Then the component $E_0\simeq \P^1$ also contains two nodes of $E$, $q_1$ and $q_2$, so that the position
of the points of attachment on $E_0$ is described by the cross-ratio
$(r_1,r_2;q_1,q_2)$. But the image of $(C,p_\bullet)$ under $\rho$ has the underlying curve $\ov{E}\cup R_1\cup R_2$
obtained by contracting $E_0$, so it does not carry the information about the relative position of $r_1,r_2,q_1,q_2$.
\end{rems}

\end{document}